\documentstyle{amsppt}
\pageheight{19cm}
\magnification=1200
\baselineskip=18pt
\nologo
\TagsOnRight
%\NoBlackBoxes
\document
\define \ts{\thinspace}
\define \oa{\frak o}
\define \g{\frak g}
\define \h{\frak h}
\define \spa{\frak {sp}}
\define \gl{\frak {gl}}
\define \sdet{\operatorname {sdet}}
\define \U{\operatorname {U}}
\define \CC{\operatorname {C}}
\define \R{\operatorname {R}}
\define \Norm{\operatorname {Norm}}
\define \Y{\operatorname {Y}}
\define \J{\operatorname {J}}
\define \Sym{{\frak S}}
\define \ZZ{{\Bbb Z}}
\define \sgn{\text{{\rm sgn}}}
\define \ot{\otimes}
\define \hra{\hookrightarrow}
\define \wh{\widehat}
\define \wt{\widetilde}
\define \Z{\operatorname {Z}}
\define \End{\operatorname {End}}
\define \tr{\text{{\rm tr\ts}}}
\define \C{\Bbb C}
\define \Proof{\noindent {\it Proof. }}

\heading{\bf A BASIS FOR REPRESENTATIONS OF}
\endheading
\heading{\bf SYMPLECTIC LIE ALGEBRAS}
\endheading 
\bigskip
\bigskip
\heading{Alexander Molev}\endheading
\bigskip
\bigskip
\bigskip
\bigskip
\noindent
Centre for Mathematics and its Applications\newline
Australian National University\newline
Canberra, ACT 0200, Australia\newline
(e-mail: molev\@pell.anu.edu.au)
\bigskip
\bigskip
\noindent
{\bf Mathematics Research Report No. MRR 012-98}
\bigskip
\bigskip
\noindent
{\bf Abstract}\newline
A basis for each finite-dimensional irreducible representation
of the symplectic Lie algebra $\spa(2n)$ is constructed.
The basis vectors are expressed in terms of the Mickelsson
lowering operators.
Explicit formulas for the matrix elements of generators of $\spa(2n)$ 
in this basis are given. The basis is
natural from the viewpoint of the representation theory of the
Yangians. The key role in the construction is played by the fact that 
the subspace of $\spa(2n-2)$-highest vectors in any 
finite-dimensional irreducible representation of $\spa(2n)$
admits a natural structure of a representation of the
Yangian $\Y(\gl(2))$.

\bigskip
\bigskip
\noindent
{\bf Mathematics Subject Classifications (1991).} 17B10, 81R10

\newpage
\noindent
{\bf 0. Introduction}
\bigskip

One of the central problems of the representation theory
is to construct a basis in the representation space and to find
the representation matrices in the basis. A solution
of this problem for the general linear Lie algebra $\gl(N)$ 
and the orthogonal Lie algebra $\oa(N)$
was given by Gelfand and Tsetlin [GT1], [GT2].
They proposed a parameterization
of basis vectors and gave formulas for the matrix
elements of the generators
of the Lie algebras in this basis. 
An explicit construction
of the Gelfand--Tsetlin basis vectors in terms of lowering
operators
was given by Zhelobenko [Z1], [Z2], Nagel--Moshinsky [NM] 
($\gl(N)$-case); Pang--Hecht [PH], Wong [Wo] ($\oa(N)$-case).
Different formulas for the lowering operators are also obtained
by Asherova--Smirnov--Tolstoy [AST2], 
Gould [G1], [G2], Nazarov--Tarasov [NT1], Molev [Mo2].

A quite different approach to construct modules over
the classical Lie algebras is developed in the papers by
King--El-Sharkaway [KS], Berele [B],
King--Welsh [KW],
Koike--Terada [KT], Proctor [P2]. It is based on the 
Weyl realization of the representations of the classical
groups in tensor spaces; see [W]. In particular,
bases in the representations of the orthogonal and symplectic
Lie algebras parameterized by $\oa(N)$-standard
or $\spa(2n)$-standard Young tableaux are constructed.
Although the subset of the standard Young tableaux is not preserved
by the action of the Lie algebra,
explicit trace relations and Garnir relations between the
Young tableaux allow one to get an algorithm for
calculation the matrix elements of the generators of the
Lie algebras.

Bases with special properties in the universal enveloping
algebra for a simple Lie algebra $\g$ and in some $\g$-modules
were constructed by Lakshmibai--Musili--Seshadri
[LMS], Littelmann [L] (monomial bases); 
De Concini--Kazhdan [CK] 
(combinatorial bases for $\text{GL}(n)$); Gelfand--Zelevinsky [GZ2],
Retakh--Zelevinsky [RZ], Mathieu [M1] (`good' bases);
Lusztig [Lu], Kashiwara [Ka] (canonical or crystal bases);
see also Mathieu [M2] for a review and more references.

The problem of constructing an analog of the Gelfand--Tsetlin
basis for the
symplectic Lie algebra $\spa(2n)$ has been addressed by many authors.
The branching rule for the reduction $\spa(2n)\downarrow\spa(2n-2)$
is obtained by Zhelobenko [Z1].
Contrary to the case of the Lie algebras $\gl(N)$ and $\oa(N)$
this reduction turns out to be not multiplicity free which
makes the problem of constructing a basis for representations of
the symplectic Lie algebra $\spa(2n)$ more complicated.

Raising and lowering
operators acting on the subspace $V(\lambda)^{+}$ of
$\spa(2n-2)$-highest vectors
in a representation $V(\lambda)$ of $\spa(2n)$ 
are constructed by Mickelsson [Mi1] (see also Bincer [Bi1], [Bi2]).
They are explicitly expressed as elements of the universal
enveloping algebra $\U(\spa(2n))$.
Applying the lowering operators 
consequently to the $\spa(2n)$-highest vector one obtains a basis
in $V(\lambda)^{+}$ and then by induction one constructs a basis
in $V(\lambda)$. However, the monomials in the
lowering operators can be chosen arbitrarily (since the operators do not
commute) and none of the bases is distinguished.
The problem of calculating
the matrix elements of generators
of $\spa(2n)$ in such a basis appears to be very difficult.

The algebra 
generated by the raising
and lowering operators, 
and more general algebras $\Z(\g,\g')$ associated with
a Lie algebra $\g$ and a reductive subalgebra $\g'\subset \g$
were studied by Mickelsson [Mi2], Van den Hombergh [Ho].
The theory of these algebras was further developed by
Zhelobenko [Z3]--[Z6] with the use of the
extremal projection method originated from [AST1]--[AST3].

A basis for the representations of the symplectic Lie algebras
was constructed by Gould and Kalnins [GK], [G3] with the use of
the restriction $\gl(2n)\downarrow \spa(2n)$.
The basis vectors are parameterized by a subset of 
the Gelfand--Tsetlin
$\gl(2n)$-patterns. Some matrix element formulas are also derived
by using the $\gl(2n)$-action.

A similar observation is made independently by Kirillov [K]
and Proctor [P1]. A description of the Gelfand--Tsetlin
patterns for $\spa(2n)$ and $\oa(N)$ can be obtained by
regarding them as fixed points of involutions of the
Gelfand--Tsetlin
patterns for the corresponding Lie algebra $\gl(N)$.

The problem of separation of multiplicities in the reduction
$\spa(2n)\downarrow \spa(2n-2)$ can be approached by 
investigating the restriction of $\spa(2n)$-modules 
to an intermediate
(non-reductive) subalgebra $\spa(2n-1)\subset\spa(2n)$. 
Such subalgebras and their representations are studied by
Gelfand--Zelevinsky [GZ1], Proctor [P1], Shtepin [S].
The separation of multiplicities can be achieved by constructing
a filtration of $\spa(2n-1)$-modules [S].

Matrix elements of generators of $\spa(2n)$ are obtained by
Wong--Yeh [WY] for certain degenerate irreducible
representations.

In this paper we give 
a construction of a weight
basis in $V(\lambda)$ and obtain
explicit formulas for the matrix elements of generators of $\g_n=\spa(2n)$
in this basis; see Theorem 1.1. Our approach is based on the
theory of Mickelsson algebras and the representation theory
of the Yangians.

It is well-known [D, Section 9.1] 
that the subspace $V(\lambda)^{+}_{\mu}$
of $\g_{n-1}$-highest vectors of a given weight $\mu$
is an irreducible representation of the centralizer algebra
$\CC_n=\U(\g_n)^{\g_{n-1}}$. However, the algebraic structure of
$\CC_n$ is very complicated which makes the problem
of studying their representations very difficult.
An approach to solve this problem is developed
by Olshanski [O3]; see also [MO].
He constructed a chain of natural homomorphisms
$$
\CC_1
\leftarrow
\CC_2
\leftarrow\cdots\leftarrow
\CC_n\leftarrow
\CC_{n+1}\leftarrow
\cdots
$$
analogous to the Harish-Chandra homomorphism [D, Section 7.4].
The projective limit of this chain is an algebra 
isomorphic to the tensor product
of an algebra of polynomials and
a quantized enveloping
algebra $\Y^-(2)$
which was called the
twisted Yangian.
(This centralizer
construction can be applied to any pairs of Lie algebras
$\frak{a}(N-M)\subset \frak{a}(N)$ of type $A$--$D$,
where $N\to\infty$ with $M$ fixed. In the result
one obtains either the Yangian $\Y(M):=\Y(\gl(M))$ 
for the Lie algebra $\gl(M)$ (see Olshanski [O1], [O2]), or
the orthogonal $\Y^+(M)$
or symplectic twisted Yangian $\Y^-(M)$; see [O3], [MO]).
In particular, one has an algebra homomorphism
$
\Y^-(2)\to\CC_n
$
so that the subspace $V(\lambda)^{+}_{\mu}$ admits a structure
of a representation of $\Y^-(2)$ which can be shown to be
irreducible.
The algebra $\Y^-(2)$ can be either defined as a subalgebra in
the Yangian $\Y(2)$ or can be presented by generators
and defining relations. The algebraic structure of the twisted Yangians
is studied in [O3] and [MNO], and their finite-dimensional irreducible
representations are described in [Mo4] in terms of the highest weights.
In particular, it is proved that
any finite-dimensional irreducible representation of $\Y^-(2)$
can be extended to the Yangian $\Y(2)$ thus providing the subspace 
$V(\lambda)^{+}_{\mu}$ with
a structure of an irreducible $\Y(2)$-module (see Theorem 5.2 below).

Lowering operators for the Yangian reduction 
$\Y(M)\downarrow \Y(M-1)$
and Gelfand-Tsetlin-type bases for representations of $\Y(M)$
were constructed in [Mo2] and [NT2] (see also [C], [NT1]).
We use a special case of these constructions
to get a Gelfand-Tsetlin-type basis 
in the $\Y(2)$-module
$V(\lambda)^{+}_{\mu}$; cf. [T], [Dr], [CP]. The basis
corresponds to an inclusion $\Y(1)\subset\Y(2)$
which can be naturally chosen by at least in two
different ways. However, to compute the action
of generators of $\g_n$ in this basis we need
to express the basis vectors in terms of the 
elements of the twisted Yangian $\Y^-(2)$.
In other words, the two inclusions
$$
\Y(1)\subset\Y(2),\qquad \Y^-(2)\subset\Y(2)
$$
must be compatible with each other in some sense 
(see Remark~4.3) which
makes the choice of the first inclusion unique and brings
the necessary rigidity into
the construction
of the basis in $V(\lambda)^{+}_{\mu}$.

To calculate the matrix elements of generators of $\g_n$ in this basis
we explicitly express the elements of the twisted Yangian $\Y^-(2)$
in terms of the
Mickelsson raising and
lowering operators. Our main instrument is Theorem 5.1 which
provides explicit formulas for the images of generators of
$\Y^-(2)$ under the natural homomorphism to the
Mickelsson algebra $\Z(\g_n,\g_{n-1})$:
$$
\Y^-(2)\to \CC_n\to \Z(\g_n,\g_{n-1}).
$$
The use of
the quadratic relations in both the algebras $\Y^-(2)$ and 
$\Z(\g_n,\g_{n-1})$ allows us to avoid long calculations.

The sections are organized as follows.
The main result (Theorem~1.1)
is formulated in Section~1. Sections~2--4
contain preliminary results which are used in the proof
of Theorem~1.1.
In Section~2 following [Z3]--[Z6] we
introduce the Mickelsson raising and lowering operators
and describe the algebraic structure of 
the algebra $\Z(\g_n,\g_{n-1})$. In
Section~3
we formulate some known results on the algebraic structure
of the Yangian $\Y(2n)$ and the twisted Yangian $\Y^-(2n)$; see 
[O3], [MNO].
In Section~4 we describe a particular case
of the construction of Gelfand--Tsetlin-type basis for
a certain class of representations of
$\Y(2)$ and $\Y^-(2)$; see [Mo2], [NT2].
Our main arguments are given in Sections~5 and 6.
We construct the highest vector and find the highest weight
for the representation $V(\lambda)^+_{\mu}$
of $\Y^-(2)$. As a corollary we obtain a proof
of the Zhelobenko branching rule for representations
of the symplectic Lie algebras [Z1]
(see also Hegerfeldt [H], King [Ki], Proctor [P2], 
Okounkov [Ok]).
In Section~6 we construct a basis in $V(\lambda)$
and derive
the formulas for the matrix elements of generators
of $\g_n$ in this basis. They have a multiplicative form
which exhibits some similarity with the Gelfand--Tsetlin
formulas in the case of $\gl(N)$ and $\oa(N)$.

\medskip

This
project was initiated in collaboration with G.~Olshanski
to whom
I would like to express my deep gratitude. I would like to thank
M.~Nazarov,
V.~Tolstoy and D.~P.~Zhelobenko for useful remarks and discussions.

\bigskip
\bigskip
\noindent
{\bf 1. Main Theorem}
\bigskip

We shall enumerate the rows and columns of $2n\times 2n$-matrices over
$\C$ by the indices $-n,-n+1,\dots,-1,1,\dots,n$. We shall also assume
throughout the paper that the index $0$ is skipped in
a summation or a product. The canonical basis $\{e_i\}$
in the space $\C^{2n}$
will be enumerated by the same set of indices.
We let the $E_{ij}$, $i,j=-n,\dots,-1,1,\dots,n$ 
denote the standard
basis of the Lie algebra $\gl(2n)$. Introduce the elements
$$
F_{ij}=E_{ij}-\theta_{ij} E_{-j,-i},\qquad 
\theta_{ij}=\sgn\ts i\cdot\sgn\ts j.
\tag 1.1
$$
The symplectic Lie algebra $\g_n=\spa(2n)$ can be identified
with the subalgebra in $\gl(2n)$
spanned by the elements $F_{ij}$, $i,j=-n,\dots,n$. They satisfy the
following symmetry property
$$
F_{-j,-i}=-\theta_{ij} F_{ij}. \tag 1.2
$$
The elements $F_{k,-k}$, $F_{-k,k}$
with $k=1,\dots,n$ and $F_{k-1,-k}$ with $k=2,\dots,n$
generate $\g_n$ as a Lie algebra.

The subalgebra $\g_{n-1}$ is spanned by the elements (1.1) with the
indices $i,j$ running over the set $\{-n+1,\dots,n-1\}$. 
Denote by $\h=\h_n$ the diagonal Cartan subalgebra in $\g_n$. 
The elements $F_{11},\dots,F_{nn}$ form a basis of $\h$. 

The finite-dimensional irreducible representations of $\g_n$
are in a one-to-one correspondence with $n$-tuples
of integers $\lambda=(\lambda_1,\dots,\lambda_n)$ 
satisfying the inequalities
$$
0\geq \lambda_1\geq\lambda_2\geq\cdots\geq \lambda_n.
$$
We denote the corresponding representation by $V(\lambda)$.
It contains a unique, up to a multiple, nonzero vector $\xi$
(the highest vector) such that
$$
\align
F_{ii}\ts\xi&=\lambda_i\ts\xi,\qquad i=1,\dots,n,\\
F_{ij}\ts\xi&=0,\qquad -n\leq i<j\leq n.
\endalign
$$
We shall sometimes 
use the numbers\footnote{The non-negative labels
$\lambda_{-n}\geq\cdots\geq\lambda_{-1}\geq 0$ 
are usually used to parameterize
the irreducible finite-dimensional representations of $\g_n$.
We have chosen to work with
positive subindices. Both parameterizations can be
easily obtained from each other.}
$\lambda_{-i}:=-\lambda_i$. They
are eigenvalues of $\xi$ with respect to the operators $F_{-i,-i}$.

The restriction of $V(\lambda)$ to the subalgebra $\g_{n-1}$
is isomorphic to a direct sum of irreducible finite-dimensional
representations $V'(\mu)$, $\mu=(\mu_1,\dots,\mu_{n-1})$
of $\g_{n-1}$ with certain
multiplicities:
$$
V(\lambda)=
\underset{\mu}\to{\bigoplus}\ts c(\mu)V'(\mu).
\tag 1.3
$$
The multiplicity $c(\mu)$ is equal to the number
of $n$-tuples of integers $(\nu_1,\dots,\nu_n)$ satisfying the
inequalities [Z1] (see also Corollary 5.3 below):
$$
\aligned
&0\geq\nu_1\geq\lambda_1\geq\nu_2\geq\lambda_2\geq \cdots\geq
\nu_{n-1}\geq\lambda_{n-1}\geq\nu_n\geq\lambda_n,\\
&0\geq\nu_1\geq\mu_1\geq\nu_2\geq\mu_2\geq \cdots\geq
\nu_{n-1}\geq\mu_{n-1}\geq\nu_n.
\endaligned
\tag 1.4
$$
Denote by $V(\lambda)^+$ the subspace of $\g_{n-1}$-highest vectors
in $V(\lambda)$:
$$
V(\lambda)^+=\{\eta\in V(\lambda)\ |\ F_{ij}\ts \eta=0,
\qquad -n<i<j<n\}.
$$
Given $\mu=(\mu_1,\dots,\mu_{n-1})$ we denote by $V(\lambda)^+_{\mu}$
the corresponding weight subspace in $V(\lambda)^+$:
$$
V(\lambda)^+_{\mu}=\{\eta\in V(\lambda)^+\ |\ F_{ii}\ts\eta=
\mu_i\ts\eta,\qquad i=1,\dots,n-1\}.
$$
We obviously have $\dim V(\lambda)^+_{\mu}=c(\mu)$.
Any nonzero vector $\eta\in V(\lambda)^+_{\mu}$ generates a
$\g_{n-1}$-submodule in $V(\lambda)$ isomorphic to $V'(\mu)$.

A parameterization of basis vectors in $V(\lambda)$ is obtained by
using its further restrictions to the subalgebras of the chain
$$
\g_1\subset\g_2\subset\cdots\subset\g_{n-1}\subset\g_n.
$$
Define the {\it pattern\/} $\Lambda$ {\it associated with\/} 
$\lambda$ as an array of integer row vectors of the form
$$
\align
&\qquad\lambda^{}_{n1}\qquad\lambda^{}_{n2}
\qquad\qquad\cdots\qquad\qquad\lambda^{}_{nn}\\
&\lambda'_{n1}\qquad \lambda'_{n2}
\qquad\qquad\cdots\qquad\qquad\lambda'_{nn}\\
&\qquad\lambda^{}_{n-1,1}\qquad\cdots\qquad\lambda^{}_{n-1,n-1}\\
&\lambda'_{n-1,1}\qquad\cdots\qquad\lambda'_{n-1,n-1}\\
&\qquad\cdots\qquad\cdots\\
&\qquad\lambda^{}_{11}\\
&\lambda'_{11}
\endalign
$$
such that the upper row coincides with $\lambda$ and 
the following inequalities hold
$$
0\geq\lambda'_{k1}\geq\lambda^{}_{k1}\geq\lambda'_{k2}\geq
\lambda^{}_{k2}\geq \cdots\geq
\lambda'_{k,k-1}\geq\lambda^{}_{k,k-1}\geq\lambda'_{kk}\geq\lambda^{}_{kk}
$$
for $k=1,\dots,n$; and
$$
0\geq\lambda'_{k1}\geq\lambda^{}_{k-1,1}\geq\lambda'_{k2}\geq
\lambda^{}_{k-1,2}\geq \cdots\geq
\lambda'_{k,k-1}\geq\lambda^{}_{k-1,k-1}\geq\lambda'_{kk}
$$
for $k=2,\dots,n$. Let us set
$$
l^{}_{ki}=\lambda^{}_{ki}-i,\qquad l'_{ki}=\lambda'_{ki}-i,\qquad
1\leq i\leq k\leq n.
\tag 1.5
$$

\bigskip
\proclaim
{\bf Theorem 1.1} There exists a basis 
$\{\zeta^{}_{\Lambda}\}$ in $V(\lambda)$ parameterized by all
patterns $\Lambda$ associated with $\lambda$ such that the action
of generators of $\g_n$ is given by the formulas
$$
\align
F_{kk}\ts \zeta^{}_{\Lambda}&=\left(2\sum_{i=1}^k\lambda'_{ki}-
\sum_{i=1}^k\lambda^{}_{ki}-\sum_{i=1}^{k-1}\lambda^{}_{k-1,i}\right)
\zeta^{}_{\Lambda},\\
F_{k,-k}\ts \zeta^{}_{\Lambda}
&=\sum_{i=1}^k A_{ki}(\Lambda)\ts
\zeta^{}_{\Lambda+\delta'_{ki}},\\
F_{-k,k}\ts \zeta^{}_{\Lambda}
&=\sum_{i=1}^k B_{ki}(\Lambda)\ts
\zeta^{}_{\Lambda-\delta'_{ki}},\\
F_{k-1,-k}\ts \zeta^{}_{\Lambda}&=\sum_{i=1}^{k-1}C_{ki}(\Lambda)\ts
\zeta^{}_{\Lambda-\delta^{}_{k-1,i}}
+\sum_{i=1}^k\sum_{j,m=1}^{k-1}D_{kijm}(\Lambda)\ts
\zeta^{}_{\Lambda+\delta'_{ki}+\delta^{}_{k-1,j}+\delta'_{k-1,m}}.
\endalign
$$
Here
$$
\align
A_{ki}(\Lambda)&=\prod_{a=1, \ts a\ne i}^k\frac{1}{l'_{ka}-l'_{ki}},\\
B_{ki}(\Lambda)&=4\ts A_{ki}(\Lambda)\ts l'_{ki}\ts 
\prod_{a=1}^k (l^{}_{ka}-l'_{ki})\prod_{a=1}^{k-1}
(l^{}_{k-1,a}-l'_{ki}),\\
C_{ki}(\Lambda)&=\frac{1}{2\ts l^{}_{k-1,i}}\prod_{a=1, \ts a\ne i}^{k-1}
\frac{1}{l_{k-1,i}^2-l_{k-1,a}^2},
\endalign
$$
and
$$
\align
&D_{kijm}(\Lambda)=A_{ki}(\Lambda)A_{k-1,m}(\Lambda)C_{kj}(\Lambda)\\
&\prod_{a=1, \ts a\ne i}^{k}(l^{}_{k-1,j}-l'_{ka})(l^{}_{k-1,j}+l'_{ka}+1)
\prod_{a=1, \ts a\ne m}^{k-1}(l^{}_{k-1,j}-l'_{k-1,a})
(l^{}_{k-1,j}+l'_{k-1,a}+1).
\endalign
$$
The arrays $\Lambda\pm\delta^{}_{ki}$ and $\Lambda\pm\delta'_{ki}$
are obtained from $\Lambda$ by replacing $\lambda^{}_{ki}$ and
$\lambda'_{ki}$
by $\lambda^{}_{ki}\pm1$ and $\lambda'_{ki}\pm1$ respectively. It is
supposed
that $\zeta^{}_{\Lambda}=0$ if the array $\Lambda$ is not a pattern.
\endproclaim

Theorem 1.1 will be proved in Sections 5 and 6.

\bigskip
\bigskip
\noindent
{\bf 2. Mickelsson algebra $\Z(\g_n,\g_{n-1})$}
\bigskip

This section contains preliminary results on 
the algebraic structure of the Mickelsson
algebra $\Z(\g_n,\g_{n-1})$; see [Z3]--[Z6] for further details.

Consider the extension 
of the universal enveloping algebra $\U(\g_n)$
$$
\U'(\g_n)=\U(\g_n)\ot_{\U(\h)} \R(\h),
$$
where $\R(\h)$ is the field of fractions of the commutative algebra
$\U(\h)$. Let $\J$
denote the left ideal in $\U'(\g_n)$ generated by 
the elements $F_{ij}$ with $-n<i<j<n$.
The {\it Mickelsson algebra\/} $\Z(\g_n,\g_{n-1})$ is the quotient algebra
of the normalizer
$$
\Norm \J=\{x\in\U'(\g_n)\ |\ \J x\subseteq \J\}
$$
modulo the two-sided ideal $\J$. 
It is an algebra over $\C$ and an $\R(\h)$-bimodule.
The algebraic structure of $\Z(\g_n,\g_{n-1})$
can be described by using the {\it extremal projection\/} $p=p_{n-1}$
for the Lie algebra $\g_{n-1}$ [AST1]--[AST3].
The projection $p$ is, up to a factor from $\R(\h_{n-1})$,
a unique element
of an extension of $\U'(\g_{n-1})$ to an algebra of formal series,
satisfying the condition
$$
F_{ij}\ts p=p\ts F_{ji}=0\qquad\text{for}\quad -n<i<j<n. \tag 2.1
$$
Explicit formulas for $p$ are given in [AST1], [Z3]. The element
$p$ is of zero weight (with respect to the adjoint action of $\h_{n-1}$)
and it can be normalized to satisfy the condition
$
p^2=p.
$
The Mickelsson algebra $\Z(\g_n,\g_{n-1})$ 
can also be defined as the image of the quotient $\U'(\g_n)/\J$
with respect to the projection $p$:
$$
\Z(\g_n,\g_{n-1})=p\left(\U'(\g_n)/\J\right).
$$
An analog of the Poincar\'e--Birkhoff--Witt theorem holds for
the algebra $\Z(\g_n,\g_{n-1})$ [Z4] so that ordered
monomials in the elements
$$
F_{n,-n},\quad F_{-n,n},\quad
pF_{in},\quad pF_{ni},\qquad i=-n+1,\dots,n-1  
$$
form a basis of
$\Z(\g_n,\g_{n-1})$ as a left or right $\R(\h)$-module.
These elements can also be given by
the following explicit formulas:
$$
\aligned
pF_{in}&=\sum_{i>i_1>\cdots>i_s>-n}
F_{ii_1}F_{i_1i_2}\cdots F_{i_{s-1}i_s}F_{i_sn}
\frac{1}{(f_i-f_{i_1})\cdots (f_i-f_{i_s})},\\
pF_{ni}&=\sum_{i<i_1<\cdots<i_s<n}
F_{i_1i}F_{i_2i_1}\cdots F_{i_si_{s-1}}F_{ni_s}
\frac{1}{(f_i-f_{i_1})\cdots (f_i-f_{i_s})}, 
\endaligned
\tag 2.2
$$
where $s=0,1,\dots$ and $f_a:=F_{aa}-a$. 
By (1.2) we have the equalities
$$
pF_{i,-n}=\sgn\ts i\cdot pF_{n,-i},\qquad pF_{-n,i}=\sgn\ts i\cdot
pF_{-i,n}.
$$
It will be
convenient to use the following normalized generators
which can be identified with elements of the universal enveloping algebra
$\U(\g_n)$: for $i=1,\dots, n-1$
$$
\aligned
z_{in}&=pF_{in}\prod_{-n<a<i}(f_i-f_a),
\qquad z_{i,-n}=pF_{i,-n}\prod_{-n<a<i}(f_i-f_a), \\
z_{ni}&=pF_{ni}\prod_{i<a<n}(f_i-f_a),
\qquad z_{-n,i}=pF_{-n,i}\prod_{i<a<n}(f_i-f_a). 
\endaligned
\tag 2.3
$$
We also set for all $i$
$$
z_{-i,n}=\sgn\ts i\cdot
z_{-n,i},\qquad z_{n,-i}=\sgn\ts i\cdot z_{i,-n}.
$$
The elements $z_{ni}$ can be written in the following
equivalent form: for $i=1,\dots,n-1$
$$
\aligned
z_{ni}&=\sum_{i<i_1<\cdots<i_s<n}(f_i-f_{j_1})\cdots (f_i-f_{j_k})\ts
F_{ni_s}F_{i_si_{s-1}}\cdots F_{i_2i_1}F_{i_1i},\\
z_{n,-i}&=\sum_{i>i_1>\cdots>i_s>-n}(f_i-f_{j_1})\cdots (f_i-f_{j_k})\ts
F_{i_s,-n}F_{i_{s-1}i_{s}}\cdots F_{i_1i_2}F_{ii_1},
\endaligned
\tag 2.4
$$
where $s=0,1,\dots$ and $\{j_1,\dots,j_k\}$ is the complementary subset
to
$\{i_1,\dots,i_s\}$ respectively in the set $\{i+1,\dots,n-1\}$ or
$\{-n+1,\dots,i-1\}$.
 
We shall use the following quadratic relations satisfied by
the elements of the algebra $\Z(\g_n,\g_{n-1})$ 
[Z4]. For all $i$
$$
\aligned
[F_{n,-n}, z_{ni}]&=0,\qquad [F_{n,-n}, z_{in}]=-2\ts z_{i,-n}, \\
[F_{-n,n}, z_{in}]&=0,\qquad [F_{-n,n}, z_{ni}]=2\ts z_{-n,i}.
\endaligned
\tag 2.5
$$
Furthermore, for $i+j\ne 0$
$$
[z_{ni},z_{nj}]=0  \tag 2.6
$$
and
$$
z_{in}\ts z_{j,-n}=z_{j,-n}\ts z_{in}\ts\frac{f_i-f_j-1}{f_i-f_j}
+z_{i,-n}\ts z_{jn}\ts\frac{1}{f_i-f_j}.
\tag 2.7
$$

Introduce the following element of $\U(\g_n)$
$$
z_{n,-n}=\sum_{n>i_1>\cdots>i_s>-n}
F_{ni_1}F_{i_1i_2}\cdots F_{i_{s-1}i_s}F_{i_s,-n}\ts
(f_{n}-f_{j_1})\cdots (f_{n}-f_{j_k}),
\tag 2.8
$$
where $s=0,1,\dots$ and $\{j_1,\dots,j_k\}$ is the complementary subset
to
$\{i_1,\dots,i_s\}$ in the set $\{-n+1,\dots,n-1\}$.
One easily checks that $z_{n,-n}$ belongs to the normalizer $\Norm \J$
and so it can be regarded as an element of $\Z(\g_n,\g_{n-1})$.
The following equivalent formula holds for $z_{n,-n}$ where
the notation of (2.8) is used:
$$
z_{n,-n}=\sum_{n>i_1>\cdots>i_s>-n}
F_{ni_1}F_{i_1i_2}\cdots F_{i_{s-1}i_s}F_{i_s,-n}\ts
(f_{-n}-f_{j_1}-1)\cdots (f_{-n}-f_{j_k}-1).
\tag 2.9
$$
To see this, we can use the standard argument 
of the extremal projection method;
see [Z3]--[Z6]. By the property (2.1)
of the extremal projection $p$ for a sequence of
indices $n>i_1>\cdots>i_s>-n$ we have 
$$
pF_{ni_1}F_{i_1i_2}\cdots F_{i_{s-1}i_s}F_{i_s,-n}
=2\ts p F_{ni_s}F_{i_s,-n},
\tag 2.10
$$
if $i_m+i_{m+1}=0$ for a certain $m$; otherwise this equals
$p F_{ni_s}F_{i_s,-n}$. This allows one to
write
the right hand side (2.8) in the form
$$
\sum_{i=1}^{n-1}pF_{ni}F_{i,-n}\ts a_i+F_{n,-n}\ts b,
\qquad a_i,b\in\R(\h).
\tag 2.11
$$
It remains to check that the coefficients $a_i$ and $b$
remain unchanged if we replace $f_n$ by $f_{-n}-1=-f_n-1$.
This can be done by
a straightforward calculation which implies that the right hand side
of (2.9) coincides with (2.11).

\bigskip
\proclaim
{\bf Proposition 2.1} We have the relations in
$\U'(\g_n)\mod \J$:
$$
\align
F_{n,-n}&=\sum_{i=-n+1}^{n}z_{ni}\ts z_{n,-i}
\prod_{\overset{\ssize a=-n+1}\to {a\ne i}}^{n}\frac{1}{f_i-f_a},
\tag 2.12\\
F_{n-1,-n}&=\sum_{i=-n+1}^{n-1}z_{n-1,i}\ts z_{i,-n}
\prod_{\overset{\ssize a=-n+1}\to {a\ne i}}^{n-1}\frac{1}{f_i-f_a},
\tag 2.13
\endalign
$$
where $z_{nn}=z_{n-1,n-1}:=1$.
\endproclaim

\Proof Both relations are proved in the same way,
so we only give a proof of (2.12).
The following equality in $\U'(\g_n)\mod \J$ is implied by
the explicit formulas for the $pF_{i,-n}$ (see (2.2)):
$$
\align
&F_{n,-n}=z_{n,-n} 
\frac{1}{(f_n-f_{-n+1})\cdots (f_n-f_{n-1})}\\
+\sum_{n>i_1>\cdots>i_s>-n}
&F_{ni_1}F_{i_1i_2}\cdots F_{i_{s-1}i_s}\cdot pF_{i_s,-n}
\frac{1}{(f_{i_s}-f_{n})(f_{i_s}-f_{i_1})
\cdots (f_{i_s}-f_{i_{s-1}})}
\endalign
$$
where $s=1,2,\dots$. Now (2.12) follows from (2.4). $\square$

\bigskip
\bigskip
\noindent
{\bf 3. Yangian $\Y(2n)$ and twisted Yangian $\Y^-(2n)$}
\bigskip

Proofs of the results formulated in this section
can be found in [MNO].

The {\it Yangian\/} $\Y(2n)=\Y(\gl(2n))$ is the
complex associative algebra with the
generators $t_{ij}^{(1)},t_{ij}^{(2)},\dots$ where 
$i,j=-n,\dots,-1,1,\dots, n$,
and the defining relations
$$
[t_{ij}(u),t_{kl}(v)]={1\over u-v}(t_{kj}(u)t_{il}(v)-t_{kj}(v)t_{il}(u)),
\tag 3.1
$$
where
$$
t_{ij}(u): = \delta_{ij} + t^{(1)}_{ij} u^{-1} + t^{(2)}_{ij}u^{-2} +
\cdots \in \Y(2n)[[u^{-1}]].
$$
One can rewrite (3.1) as a
{\it ternary relation\/}
for the matrix
$$
T(u):=\sum_{i,j} t_{ij}(u)\ot E_{ij}
\in \Y(2n)[[u^{-1}]]\ot \End\C^{2n}.
$$
To do this introduce the following notation.
For an operator $X\in \End\C^{2n}$ and a number $m=1,2,\ldots$ we set
$$
X_k:=1^{\otimes (k-1)}\otimes X\otimes 1^{\otimes (m-k)}
\in\left(\End\C^{2n}\right)^{\otimes m}, \quad 1\leq k\leq m. \tag 3.2
$$
If $X\in(\End\C^{2n})^{\otimes 2}$ then for any $k, l$\  such that 
$1\leq k, l\leq m$\ and $k\ne l$, we denote by $X_{kl}$ the operator in 
$(\C^{2n})^{\otimes m}$
which acts as $X$ in the product of $k$th and $l$th copies and as
1 in all other copies. That is, 
$$
X=\sum_{r,s,t,u} a_{rstu} E_{rs}\otimes E_{tu}\quad
\Rightarrow \quad X_{kl}=\sum_{r,s,t,u} a_{rstu} (E_{rs})_k\ (E_{tu})_l, 
\tag 3.3
$$
where $a_{rstu}\in\C$.
The ternary relation has the form
$$
R(u-v)T_1(u)T_2 (v) = T_2(v)T_1(u)R(u-v), \tag 3.4
$$
where $R(u)=R_{12}(u)=1-u^{-1}P$ and $P$ is the permutation operator in
$\C^{2n}\ot \C^{2n}$.

The Yangian $\Y(2n)$ is a Hopf algebra with the 
coproduct 
$$
\Delta (t_{ij}(u))=\sum_{a=-n}^n t_{ia}(u)\ot
t_{aj}(u).\tag 3.5
$$

The {\it twisted Yangian\/} $\Y^{-}(2n)$
corresponding to the
symplectic Lie algebra $\g_n=\spa(2n)$ is defined as follows.
By $X\mapsto X^t$ we will denote the matrix 
transposition such that
$
(E_{ij})^t=\theta_{ij}E_{-j,-i}.
$
Introduce the matrix $S(u)=(s_{ij}(u))$ by setting
$
S(u):=T(u)T^t(-u),
$
or, in terms of matrix elements,
$$
s_{ij}(u)=\sum_{a=-n}^n \theta_{aj}t_{ia}(u)t_{-j,-a}(-u). \tag 3.6
$$
Write
$
s_{ij}(u)=\delta_{ij}+s_{ij}^{(1)}u^{-1}+s_{ij}^{(2)}u^{-2}+\cdots.
$
The twisted Yangian $\Y^{-}(2n)$ is the subalgebra of $\Y(2n)$ generated by
the elements $s_{ij}^{(1)},s_{ij}^{(2)},\dots$, where $-n\leq i,j\leq n$.

The matrix $S(u)$ satisfies the following
{\it quaternary relation\/} and {\it symmetry relation\/} which
follow from (3.4):
$$
\aligned
R(u-v)S_1(u)R^t(-u-v)S_2(v)={}&S_2(v)R^t(-u-v)S_1(u)R(u-v), \\
S^t(-u)={}&\frac{2u-1}{2u}\ts S(u)+\frac{1}{2u}\ts S(-u).
\endaligned
\tag 3.7
$$
Here we use the notation (3.2),
where $R^t(u)$ is obtained from $R(u)$ by applying the transposition $t$
in either of the two copies of $\End\C^{2n}$:
$$
R^t(u)=1-u^{-1}\sum_{i,j}\theta_{ij}E_{-j,-i}\ot E_{ji}.
$$
Relations (3.7) are defining relations for the algebra $\Y^{-}(2n)$ and
they can be rewritten in terms of 
the generating series $s_{ij}(u)$ as follows:
$$
\aligned
[s_{ij}(u),s_{kl}(v)]={}&{1\over
u-v}(s_{kj}(u)s_{il}(v)-s_{kj}(v)s_{il}(u))\\
-{}&{1\over u+v}(\theta_{k,-j}s_{i,-k}(u)s_{-j,l}(v)-
\theta_{i,-l}s_{k,-i}(v)s_{-l,j}(u))\\
+{}&{1\over u^2-v^2}(\theta_{i,-j}s_{k,-i}(u)s_{-j,l}(v)-
\theta_{i,-j}s_{k,-i}(v)s_{-j,l}(u)) 
\endaligned
\tag 3.8
$$
and
$$
\theta_{ij}s_{-j,-i}(-u)=\frac{2u-1}{2u}\ts
s_{ij}(u)+\frac{1}{2u}\ts s_{ij}(-u). 
\tag 3.9
$$
This allows one to regard $\Y^-(2n)$ as an abstract algebra with
generators $s_{ij}^{(r)}$ and the relations (3.8), (3.9).

The mapping
$
F_{ij}\mapsto s_{ij}^{(1)}
$
defines an inclusion $\U(\g_n)\hra\Y^{-}(2n)$ while
the mapping
$$
s_{ij}(u)\mapsto \delta_{ij}+\frac{F_{ij}}{u-1/2}
\tag 3.10
$$
defines an algebra homomorphism
$
\Y^{-}(2n)\to \U(\g_n).
$
Any even formal series $c(u)\in 1+u^{-2}\C[[u^{-2}]]$ defines
an automorphism of $\Y^{-}(2n)$ given by
$$
s_{ij}(u)\mapsto c(u)\ts s_{ij}(u).
\tag 3.11
$$

The
{\it Sklyanin determinant\/}
$\sdet S(u)$
is a formal series in $u^{-1}$ with
coefficients from the center of the algebra
$\Y^{-}(2n)$. It can be defined by the formula (see (3.2), (3.3)):
$$
\align
A_{2n}S_1(u)R^t_{12}&\cdots R^t_{1,2n}S_2(u-1)
R^t_{23}\cdots R^t_{2,2n}S_3(u-2)\\
{}&\cdots S_{2n-1}(u-2n+2)R^t_{2n-1,2n}S_{2n}(u-2n+1)=\sdet S(u)A_{2n},
\tag 3.12
\endalign
$$
where
$R^t_{ij}:=R^t_{ij}(-2u+i+j-2)$, and $A_{2n}$ is the normalized
antisymmetrizer
in the tensor space $(\C^{2n})^{\otimes {2n}}$ so that $A_{2n}^2=A_{2n}$.
Explicit formulas for
$\sdet S(u)$ are given in [O3], [MNO], [Mo3].

The {\it Sklyanin comatrix\/}
$\widehat S(u)=(\widehat s_{ij}(u))$ is defined by
$$
\sdet S(u)=\widehat S(u)S(u-2n+1).\tag 3.13
$$
The mapping
$$
S(u)\mapsto \frac{2u+1}{2u-2n+1}\ts \wh S(-u+n-1)
\tag 3.14
$$
defines an automorphism of the algebra $\Y^{-}(2n)$; 
see [Mo4, Proposition 1.1].

Let us denote by $S^{(n-1)}(u)$
and $\wt S(u)$ the submatrices
of $S(u)$
whose rows and columns are enumerated by the sets of indices
$\{-n+1,\dots,-1,1,\dots,n-1\}$ and 
$\{-n+1,\dots,-1,1,\dots,n\}$ respectively. Introduce
the $nn$th {\it quasi-determinant\/}
of the matrix $\wt S(u)$ by
$$
|\wt S(u)|_{nn}=\left(\left(\wt S(u)^{-1}\right)_{nn}\right)^{-1};
$$
see [GKLLRT].
We shall need the following expression for the matrix
element $\wh s_{nn}(u)$ of the Sklyanin comatrix $\wh S(u)$.

\bigskip
\proclaim
{\bf Proposition 3.1} We have the formula
$$
\wh s_{nn}(u)=\frac{2u+1}{2u-1}\ts
|\wt S(-u)|_{nn}\ts\sdet S^{(n-1)}(u-1).
\tag 3.15
$$
\endproclaim

\Proof Multiplying both sides of (3.12) by $S_{2n}^{-1}(u-2n+1)$ 
from the right and using (3.13)
we obtain the relation
$$
\align
A_{2n}S_1(u)R^t_{12}\cdots R^t_{1,2n}S_2(u-1)&
R^t_{23}\cdots R^t_{2,2n}S_3(u-2)\\
{}\cdots \ts &S_{2n-1}(u-2n+2)R^t_{2n-1,2n}
=A_{2n}\wh S_{2n}(u).\tag 3.16
\endalign
$$
It can be easily verified by
using the symmetry relation (3.9) (see also [Mo3]) that
$$
A_{2n}S_1(u)R^t_{12}\cdots R^t_{1,2n}=
\frac{2u+1}{2u-1}\ts A_{2n}S_1^t(-u). 
$$
Denote by $A_{2n}^{(2)}$ the normalized
antisymmetrizer corresponding to
the subgroup $\Sym_{\{2,\dots,2n\}}$ of the symmetric group
$\Sym_{2n}$.
Clearly, $A_{2n}=A_{2n}A_{2n}^{(2)}$. Note
that $A_{2n}^{(2)}$ is permutable with $S_1^t(-u)$, while $R^t_{ij}$
is permutable with
$R^t_{kl}$ and $S_k(u)$ provided that the indices $i,j,k,l$ are distinct.
So, we can rewrite formula (3.16) in the form:
$$
\align
\frac{2u+1}{2u-1}\ts A_{2n}S_1^t(-u)A_{2n}^{(2)}&S_2(u-1)
R^t_{23}\cdots R^t_{2,{2n}-1}S_3(u-2)\\
\cdots \ts &S_{{2n}-1}(u-{2n}+2)R^t_{2,{2n}}\cdots R^t_{{2n}-1,{2n}}
=A_{2n}\wh S_{2n}(u).\tag 3.17
\endalign
$$
Let us apply the operators in
both sides of this formula to the vector
$v_i=e_{-i}\ot e_{-n+1}\ot e_{-n+2}\ot\cdots\ot e_{n-1}\ot e_n$, where
$i\in\{-n+1,\dots,n\}$. For the right hand side we clearly obtain
$$
A_{2n}\wh S_{2n}(u)v_i=
\delta_{in}\ts\wh s_{nn}(u)\ts\zeta,
\tag 3.18
$$
where $\zeta:=A_{2n}(e_{-n}\ot e_{-n+1}
\ot\cdots\ot e_n)$.
To calculate the left hand side we note first that
$$
R^t_{2,{2n}}\cdots R^t_{{2n}-1,{2n}}v_i=v_i.
$$
Further, let us introduce the formal series
$$
\Phi_{a_2,\dots,a_{2n-1}}(u-1)\in\Y^{-}(2n)[[u^{-1}]],\quad
-n\leq a_i\leq n,
$$
as follows:
$$
\align
A_{2n}^{(2)}S_2(u-1)
R^t_{23}\cdots R^t_{2,{2n}-1}S_3(u-2)
\cdots S_{{2n}-1}(u-{2n}+2)&(e_{-n+1}\ot\cdots\ot e_{n-1})\\
{}=\sum_{a_2,\dots,a_{{2n}-1}}\Phi_{a_2,\dots,a_{{2n}-1}}(u-1)
&(e_{a_2}\ot\cdots\ot e_{a_{{2n}-1}}).
\endalign
$$
In particular,
$$
({2n}-2)!\ \Phi_{-n+1,\dots,n-1}(u-1)=\sdet S^{(n-1)}(u-1),\tag 3.19
$$
and the series $\Phi_{a_2,\dots,a_{{2n}-1}}(u-1)$
is skew symmetric with respect to permutations of the indices
$a_2,\dots,a_{2n-1}$; see [MNO, Section 4]. 
This allows us to write the left hand side of
(3.17) applied to $v_i$ in the form:
$$
\align
\frac{2u+1}{2u-1}\ts&({2n}-2)!\ts\sum_{k=1}^{{2n}-1}(-1)^{k-1}\ts
s_{b_k,-i}^t(-u)\ts
\Phi_{b_1,\dots,\wh {b}_k,\dots,b_{{2n}-1}}(u-1)\ts\zeta\\
=\frac{2u+1}{2u-1}\ts&({2n}-2)!\ts\theta_{in}\sum_{k=1}^{{2n}-1}
s_{i,-b_k}(-u)\ts(-1)^{k-1}\ts
\theta_{-b_k,n}\ts\Phi_{b_1,\dots,\wh {b}_k,\dots,b_{{2n}-1}}(u-1)\ts\zeta,
\endalign
$$
where $(b_1,\dots,b_{{2n}-1})=(-n,-n+1,\dots,n-1)$ and the hat
indicates the index to be omitted.
Put
$$
\Phi_{-b_k}(u-1):=({2n}-2)!\ts(-1)^{k-1}\ts
\theta_{-b_k,n}\ts\Phi_{b_1,\dots,\wh {b}_k,\dots,b_{{2n}-1}}(u-1).
$$
Then, taking into account (3.18), we get the following matrix
relation:
$$
\frac{2u+1}{2u-1}\ts\wt S(-u)
\left( \matrix\Phi_{-n+1}(u-1)\\
\vdots\\
\Phi_n(u-1)
\endmatrix\right)=
\left( \matrix 0\\
\vdots\\
\wh s_{nn}(u)
\endmatrix\right)
$$
Multiplying its both sides by the matrix $\wt S(-u)^{-1}$ from the
left and comparing the $n$th coordinates of the vectors,
we obtain using (3.19) that
$$
\frac{2u+1}{2u-1}\ts\sdet S^{(n-1)}(u-1)=
\bigl(\wt S(-u)^{-1}\bigr)_{nn}\ts\wh s_{nn}(u),
$$
which implies (3.15). $\square$

\bigskip
\bigskip
\noindent
{\bf 4. Representations of the
algebras $\Y(2)$ and $\Y^-(2)$}
\bigskip

Here we formulate some necessary results on
representations of the algebras $\Y(2)$ and $\Y^-(2)$;
see [T], [Dr], [CP], [NT2], [Mo2], [Mo4].
Having in mind their applications in Sections~5 and 6 we shall
enumerate the generators of $\Y(2)$ and $\Y^-(2)$,
as well as rows and columns of $2\times 2$-matrices, by the symbols
$-n,n$ instead of the usual $-1,1$.

A representation of the Yangian $\Y(2)$ is called
{\it highest weight\/} if it is
generated by a nonzero vector $\eta$ such that
$$
\align
t_{i,i}(u)\ts\eta&=\lambda_i(u)\ts\eta,\qquad i=-n,n,\\
t_{-n,n}(u)\ts\eta&=0,
\endalign
$$
for certain formal series $\lambda_i(u)\in1+u^{-1}\C[[u^{-1}]]$.
The pair $(\lambda_{-n}(u),\lambda_n(u))$ is called the 
{\it highest weight\/} of the representation.
Given arbitrary series $\lambda_{-n}(u),\lambda_n(u)$ there exists
a unique, up to an isomorphism, irreducible highest weight representation
of
$\Y(2)$ with the highest weight $(\lambda_{-n}(u),\lambda_n(u))$
which will be denoted by $L(\lambda_{-n}(u),\lambda_n(u))$.

Similarly, 
a representation of the Yangian $\Y^-(2)$ is {\it highest weight\/} if it
is
generated by a nonzero vector $\eta$ such that
$$
\align
s_{n,n}(u)\ts\eta&=\mu(u)\ts\eta,\\
s_{-n,n}(u)\ts\eta&=0,
\endalign
$$
for a certain formal series $\mu(u)\in1+u^{-1}\C[[u^{-1}]]$
called the {\it highest weight\/} of the representation.
Given an arbitrary series $\mu(u)$ there exists
a unique, up to an isomorphism, irreducible highest weight representation
of
$\Y^-(2)$ with the highest weight $\mu(u)$
which will be denoted by $V(\mu(u))$. Every irreducible
finite-dimensional representation of the algebra $\Y^-(2)$
is isomorphic to a unique $V(\mu(u))$.

Given a pair of complex numbers $(\alpha,\beta)$ 
such that $\alpha-\beta\in\ZZ_+$
we denote by 
$L(\alpha,\beta)$ the irreducible representation of the Lie algebra
$\gl(2)$ with the highest weight $(\alpha,\beta)$ with respect to the
upper triangular Borel subalgebra. We have 
$\dim L(\alpha,\beta)=\alpha-\beta+1$. We may regard $L(\alpha,\beta)$ as
a $\Y(2)$-module by using the algebra homomorphism $\Y(2)\to\U(\gl(2))$
given by
$$
t_{ij}(u)\mapsto \delta_{ij}+E_{ij}u^{-1},\qquad i,j\in\{-n,n\}.
\tag 4.1
$$
The coproduct (3.5) allows one to construct representations of
$\Y(2)$ of the form
$$
L=L(\alpha_1,\beta_1)\ot\cdots\ot L(\alpha_n,\beta_n).
\tag 4.2
$$
One easily obtains from (3.5) and (4.1) that the tensor product
$$
\eta=\omega_1\ot\cdots\ot\omega_n \tag 4.3
$$ 
of the highest weight vectors 
$\omega_i$ of the $L(\alpha_i,\beta_i)$ generates a highest weight
submodule in $L$ with the highest weight $(\alpha(u),\beta(u))$,
where
$$
\aligned
\alpha(u)&=(1+\alpha_1u^{-1})\cdots (1+\alpha_nu^{-1}),\\
\beta(u)&=(1+\beta_1u^{-1})\cdots (1+\beta_nu^{-1}).
\endaligned
\tag 4.4
$$
Hence, if the representation $L$ is irreducible then it is
isomorphic to $L(\alpha(u),\beta(u))$. A criterion of irreducibility of
representation (4.2) is given by Chari and Pressley [CP] and can be also
deduced from results of Tarasov [T] (see [Mo4]). To formulate
the result, with each
$L(\alpha,\beta)$ associate the {\it string\/}
$$
S(\alpha,\beta)=\{\beta,\beta+1,\dots,\alpha-1\}\subset\C.
$$
We say that two strings $S_1$ and $S_2$ are 
{\it in general position\/}
if either

(i) $S_1\cup S_2$ is not a string, or

(ii) $S_1\subseteq S_2$, or $S_2\subseteq S_1$.

The representation (4.2) of $\Y(2)$ is irreducible 
if and only if the strings
$S(\alpha_i,\beta_i)$ and $S(\alpha_j,\beta_j)$ are in
general position for all $i<j$ [CP].

The tensor product (4.2) can also be regarded as a representation
of the subalgebra $\Y^-(2)\subset\Y(2)$. The following criterion of its
irreducibility is given in [Mo4] and will be used in Section~5.

\bigskip
\proclaim
{\bf Proposition 4.1} The representation (4.2) of $\Y^-(2)$
is irreducible if and only if each pair of strings
$(S(\alpha_i,\beta_i), S(\alpha_j,\beta_j))$ and
$(S(\alpha_i,\beta_i), S(-\beta_j,-\alpha_j))$ is in general position
for all $i<j$. $\square$
\endproclaim

If the representation $L$ of $\Y^-(2)$
defined in (4.2) is irreducible then
by (3.6) it is isomorphic to $V(\mu(u))$ with
$$
\mu(u)=\alpha(-u)\beta(u), \tag 4.5
$$
where $\alpha(u)$ and $\beta(u)$ are given by (4.4).

It follows from (3.5) and (4.1) that the elements $t_{ij}^{(r)}\in\Y(2)$
with $r\geq n+1$ act as zero operators in (4.2). Therefore,
the operators
$$
T_{ij}(u)=u^n\ts t_{ij}(u) \tag 4.6
$$
are polynomials in $u$. By (3.6) the same is true for the
operators
$
S_{ij}(u)=u^{2n}\ts s_{ij}(u).
$
Note that the defining relations (3.1) allow us
to rewrite the formula (3.6) for $s_{n,-n}(u)$ in the form
$$
s_{n,-n}(u)=\frac{u+1/2}{u}\left(t_{n,-n}(u)t_{nn}(-u)-
t_{n,-n}(-u)t_{nn}(u)\right).
$$
Therefore
we may introduce another polynomial operator in $L$ by
$$
S^{\natural}_{n,-n}(u)=\frac{1}{u+1/2}\ts S_{n,-n}(u)
=\frac{(-1)^n}{u}\left(T_{n,-n}(u)T_{nn}(-u)-
T_{n,-n}(-u)T_{nn}(u)\right). \tag 4.7
$$
Note that by (3.8) we have 
$[S^{\natural}_{n,-n}(u), S^{\natural}_{n,-n}(v)]=0$. 

Let $\gamma=(\gamma_1,\dots,\gamma_n)$ be an $n$-tuple of complex numbers
such that
$$
\alpha_i-\gamma_i\in\ZZ_+,\qquad \gamma_i-\beta_i\in\ZZ_+,\qquad
i=1,\dots,n.
\tag 4.8
$$
Introduce the following vector in $L$
$$
\eta_{\gamma}=\prod_{i=1}^n S^{\natural}_{n,-n}(-\gamma_i+1)\cdots 
S^{\natural}_{n,-n}(-\beta_i-1)S^{\natural}_{n,-n}(-\beta_i)\ts\eta,
\tag 4.9
$$
where $\eta$ is defined in (4.3).

\bigskip
\proclaim
{\bf Proposition 4.2} Suppose that the representation $L$
of $\Y^-(2)$
given by (4.2) is irreducible
and the strings $S(\alpha_i,\beta_i)$ satisfy the condition
$$
S(\alpha_i,\beta_i)\cap S(\alpha_j,\beta_j)=\emptyset
\qquad\text{for}\quad i\ne j. \tag 4.10
$$
Then the vectors $\eta_{\gamma}$ with $\gamma$ satisfying (4.8)
form a basis in $L$. Moreover, one has the formulas
$$
\align
T_{nn}(u)\ts\eta_{\gamma}&=(u+\gamma_1)\cdots (u+\gamma_n)
\ts\eta_{\gamma}, \tag 4.11\\
T_{n,-n}(-\gamma_i)\ts \eta_{\gamma}&=
\frac12\prod_{a=1,\ts a\ne i}^n\frac{1}{-\gamma_i-\gamma_a}\ts
\eta_{\gamma+\delta_i},\tag 4.12\\
T_{-n,n}(-\gamma_i)\ts \eta_{\gamma}&=
-2\ts\prod_{k=1}^n(\alpha_k-\gamma_i+1)(\beta_k-\gamma_i)\cdot
\prod_{a=1,\ts a\ne i}^n (-\gamma_i-\gamma_a+1)
\ts\eta_{\gamma-\delta_i}, \tag 4.13
\endalign
$$
where $\delta_i$ is the $n$-tuple which has $1$ on the $i$th position
and zeroes as remaining entries; it is assumed that $\eta_{\gamma}=0$
if $\gamma$ does not satisfy (4.8).
\endproclaim

\Proof Since $L$ is irreducible as a $\Y(2)$-module
it is isomorphic to the highest weight
representation $L(\alpha(u),\beta(u))$.
For each $\gamma$ satisfying (4.8) introduce 
the vector
$$
\wt{\eta}_{\gamma}=\prod_{i=1}^n T_{n,-n}(-\gamma_i+1)\cdots 
T_{n,-n}(-\beta_i-1)T_{n,-n}(-\beta_i)\ts\eta.
$$
The vectors $\{\wt{\eta}_{\gamma}\}$
form a basis in $L$ and the following relations hold:
$$
\align
T_{nn}(u)\ts\wt{\eta}_{\gamma}&=(u+\gamma_1)\cdots (u+\gamma_n)
\ts\wt{\eta}_{\gamma}, \tag 4.14\\
T_{n,-n}(-\gamma_i)\ts 
\wt{\eta}_{\gamma}&=\wt{\eta}_{\gamma+\delta_i},\tag 4.15\\
T_{-n,n}(-\gamma_i)\ts \wt{\eta}_{\gamma}&=
-\prod_{k=1}^n(\alpha_k-\gamma_i+1)(\beta_k-\gamma_i)
\ts\wt{\eta}_{\gamma-\delta_i}. \tag 4.16
\endalign
$$
This is a special case of the construction of
Gelfand--Tsetlin-type bases for representations of the Yangian
$\Y(m)$ [Mo2], [NT2] (see also [T], [NT1]).

The formulas (4.7), (4.14) and (4.15) imply that
$$
S^{\natural}_{n,-n}(-\gamma_i)\ts\wt{\eta}_{\gamma}=2\ts 
\prod_{a=1,\ts a\ne i}^n (-\gamma_i-\gamma_a) \ts 
\wt{\eta}_{\gamma+\delta_i}.
\tag 4.17
$$
Hence, for each $\gamma$ the vectors $\eta_{\gamma}$ and
$\wt{\eta}_{\gamma}$ coincide up to a nonzero multiple.
This proves (4.11). Now (4.12) and (4.13) follow
from (4.15) and (4.16).
$\square$
\medskip

\noindent
{\it Remark 4.3.\ } The above proof of (4.17) relies on the fact that the 
$\wt{\eta}_{\gamma}$ are eigenvectors for the operators $T_{nn}(u)$;
see (4.7). That is, the Gelfand--Tsetlin-type basis 
$\{\wt{\eta}_{\gamma}\}$ in $L$ corresponds to the inclusion
$\Y(1)\subset \Y(2)$ with $\Y(1)$ generated by the coefficients
of $t_{nn}(u)$.

\bigskip
\bigskip
\noindent
{\bf 5. Yangian action on $V(\lambda)^+_{\mu}$}
\bigskip

Let us introduce the
$2n\times 2n$-matrix $F=(F_{ij})$ whose $ij$th entry is
the element $F_{ij}\in\g_n$  (see (1.1)) and set
$$
\Cal F(u)=1+\frac{F}{u-1/2}.
$$
Denote by $\wh {\Cal F}(u)$ the image of the Sklyanin comatrix
$\wh S(u)$ under the homomorphism $S(u)\mapsto \Cal F(u)$; see (3.10).
By (3.11) and (3.14) the mapping
$$
\pi:
S(u)\mapsto c(u)\ts \frac{2u+1}{2u-2n+1}\ts \wh {\Cal F}(-u+n-1),
$$
where
$$
c(u)=\prod_{k=1}^n (1-(k-1/2)^2\ts u^{-2}),
$$
defines a homomorphism $\Y^-(2n)\to \U(\g_n)$; cf. [O3], [MO].
The series
$$
(1+nu^{-1})^{-1}\sdet S(u+n-1/2)
$$
is even in $u$ (see [MNO, Section 4.11]) and so by (3.14)
the image of the generator
$s_{ij}^{(1)}$ with respect to $\pi$ coincides with $F_{ij}$.
By (3.8) we then have
$$
[F_{ij},s^{\pi}_{kl}(u)]=\delta_{kj}s^{\pi}_{il}(u)-
\delta_{il}s^{\pi}_{kj}(u)
-\theta_{k,-j}\delta_{i,-k}s^{\pi}_{-j,l}(u)+
\theta_{i,-l}\delta_{-l,j}s^{\pi}_{k,-i}(u),
\tag 5.1
$$
where $s^{\pi}_{ij}(u):=\pi(s_{ij}(u))$. This implies that
the image of the restriction of $\pi$ to the
subalgebra $\Y^-(2)$ generated by the elements $s_{ij}(u)$
with $i,j\in\{-n,n\}$ is contained in the centralizer 
$\CC_n=\U(\g_n)^{\g_{n-1}}$ and thus defines an algebra homomorphism
$$
\pi:\Y^-(2)\to \CC_n. \tag 5.2
$$
However,
the subspace $V(\lambda)^+_{\mu}$ is an irreducible representation
of the centralizer $\CC_n$; see [D, Section 9.1].
It follows from [MO, Proposition 4.9] that
the algebra $\CC_n$ is generated by the image of $\pi$ and
the center of $\U(\g_n)$. Since the elements of the center
of $\U(\g_n)$ act as scalar operators in $V(\lambda)$,
the $\Y^-(2)$-module $V(\lambda)^+_{\mu}$ defined via 
the homomorphism (5.2) is irreducible.

Note that $\CC_n$ is a subalgebra in the 
normalizer $\Norm \J$ (see Section 2):
$$
\CC_n\hra \Norm \J.
$$
Thus, using (5.2) and the definition of the
Mickelsson algebra $\Z(\g_n,\g_{n-1})$ we
obtain an algebra homomorphism which we still denote by $\pi$:
$$
\pi:\Y^-(2)\to \Z(\g_n,\g_{n-1}). \tag 5.3
$$
In other words, the elements of $\Y^-(2)$, as operators in
the space $V(\lambda)^+_{\mu}$, can be expressed as elements
of the Mickelsson algebra $\Z(\g_n,\g_{n-1})$. An explicit form of the
images of the generators of $\Y^-(2)$ under
the homomorphism (5.3)
is given in the following theorem.

Introduce the polynomials 
$Z_{ij}(u)$,\ \  $i,j\in\{-n,n\}$
with coefficients in the 
Mickelsson algebra
$\Z(\g_n,\g_{n-1})$:
$$
\align
Z_{n,-n}(u)&=
\sum_{i=1}^{n-1}z_{ni}\ts z_{n,-i}\prod_{a=1,\ts a\ne i}^{n-1}
\frac{u^2-g_a^2}{g_i^2-g_a^2}+
F_{n,-n}\ts 
\prod_{a=1}^{n-1}(u^2-g_a^2),
\tag 5.4\\
Z_{-n,n}(u)&=
\sum_{i=1}^{n-1}z_{-i,n}\ts z_{in}\prod_{a=1,\ts a\ne i}^{n-1}
\frac{u^2-g_a^2}{g_i^2-g_a^2}+
F_{-n,n}\ts
\prod_{a=1}^{n-1}(u^2-g_a^2),
\tag 5.5\\
Z_{n,n}(u)&=
\sum_{i=-n+1}^{n-1}z_{ni}\ts z_{-n,-i}
\prod_{\overset{\ssize a=-n+1}\to {a\ne i}}^{n-1}
\frac{u+g_a}{g_i-g_a}+(u+g_n)
\prod_{a=-n+1}^{n-1}(u+g_a),
\tag 5.6\\
Z_{-n,-n}(u)&=-\sum_{i=-n+1}^{n-1}z_{-i,n}\ts z_{i,-n}
\prod_{\overset{\ssize a=-n+1}\to {a\ne i}}^{n-1}
\frac{u+g_a}{g_i-g_a}+(u+g'_n)
\prod_{a=-n+1}^{n-1}(u+g_a),
\tag 5.7
\endalign
$$
where $g_i:=f_i+1/2=F_{ii}-i+1/2$ and $g'_n=-g_n-2n+1$.

\bigskip
\proclaim
{\bf Theorem 5.1} The images of the generators $s_{ij}(u)$,
$i,j\in\{-n,n\}$
of $\ts\Y^-(2)$
under the homomorphism (5.3) are given by
the formulas
$$
\alignedat2
s_{-n,-n}(u)&\mapsto\frac{u+1/2}{u^{2n}}\ts Z_{-n,-n}(u),\qquad
s_{-n,n}(u)&&\mapsto\frac{u+1/2}{u^{2n}}\ts Z_{-n,n}(u),\\
s_{n,-n}(u)&\mapsto\frac{u+1/2}{u^{2n}}\ts Z_{n,-n}(u),\qquad\quad
s_{n,n}(u)&&\mapsto\frac{u+1/2}{u^{2n}}\ts Z_{n,n}(u).
\endalignedat
\tag 5.8
$$
\endproclaim

\Proof Consider first the generator $s_{n,n}(u)$. 
Using Proposition 3.1 we obtain the following formula
for the $nn$th matrix element of the matrix $\wh{\Cal F}(u-1/2)$:
$$
\wh{\Cal F}(u-1/2)_{nn}=\frac{u}{u-1}\ts |1-\wt F\ts u^{-1}|_{nn}\ts
\sdet \Cal F^{(n-1)}(u-3/2), \tag 5.9
$$
where $\wt F$ is the submatrix of $F$ obtained by 
removing the row and column enumerated by $-n$, and
$\sdet \Cal F^{(n-1)}(u)$ is the image of the Sklyanin determinant
$\sdet S^{(n-1)}(u)$ under the homomorphism (3.10).
Using the combinatorial interpretation of the
quasi-determinant $|1-\wt F\ts u^{-1}|_{nn}$
[GKLLRT, Proposition 7.20]
we obtain the formula:
$$
|1-\wt F\ts u^{-1}|_{nn}=1-\sum_{k=1}^{\infty}F^{(k)}_{nn}\ts u^{-k},
\tag 5.10
$$
where
$$
F^{(k)}_{ab}=\sum F^{}_{ai_1}F^{}_{i_1i_2}\cdots F^{}_{i_{k-1}b},
$$
summed over all values of the indices 
$i_m\in\{-n+1,\dots,-1,1,\dots,n-1\}$. Let us show that
for $k\geq 2$ we have the equality in $\Z(\g_n,\g_{n-1})$:
$$
F^{(k)}_{nn}=\sum_{i=-n+1}^{n-1}p\ts F^{(k-1)}_{ni}\cdot
\prod_{i<a<n}(f_i-f_a)
\cdot pF^{}_{in}\ts\prod_{i<a<n}\frac{1}{f_i-f_a},
\tag 5.11
$$
where $p$ is the extremal projection for $\g_{n-1}$; see Section 2.
First, we note that
$$
[F^{(k-1)}_{ni}, F^{}_{ab}]=\delta_{ai}F^{(k-1)}_{nb}-
\theta_{ab}\ts\delta_{i,-b}F^{(k-1)}_{n,-a}, \tag 5.12
$$
where $i,a,b\in\{-n+1,\dots,n-1\}$.
Indeed, the coefficients of $\sdet S(u)$ are central in $\Y^-(2)$
and so by (3.13) and the definition of $\pi$
the relations (5.1) hold for the $s^{\pi}_{ij}(u)$ replaced with
$(1-Fu^{-1})_{ij}^{-1}$. Taking the
coefficient at $u^{-m}$ we get
the well-known formula
for the commutator $[F^{}_{ij},(F^m_{})^{}_{kl}]$ which implies (5.12).
Now we transform the right hand side of (5.11) by using
the explicit formula
for the generators $pF^{}_{in}$ given in (2.2). By the property (2.1)
of $p$ we obtain from (5.12) that
if $i>i_1>\cdots>i_s>-n$ then
$$
pF^{(k-1)}_{ni}\cdot F^{}_{ii_1}F^{}_{i_1i_2}
\cdots F^{}_{i_{s-1}i_s}F^{}_{i_sn}
=2\ts p F^{(k-1)}_{ni_s}F^{}_{i_sn},
$$
if $i+i_1=0$ or $i_m+i_{m+1}=0$ for a certain $m$; 
otherwise this equals
$p F^{(k-1)}_{ni_s}F^{}_{i_sn}$ (cf. (2.10)).
Using this, we verify by a straightforward calculation
that the coefficient
at each product $pF^{(k-1)}_{ni} F^{}_{in}$
on the right hand side of (5.11) equals $1$,
and so it is given by
$$
\sum_{i=-n+1}^{n-1}pF^{(k-1)}_{ni} F^{}_{in}=pF^{(k)}_{nn}=F^{(k)}_{nn},
$$
which proves (5.11).
Further, we have in $\Z(\g_n,\g_{n-1})$
$$
pF^{(k-1)}_{ni}=pF^{}_{ni} (f_i+n-1)^{k-2},\qquad i=-n+1,\dots,n-1.
\tag 5.13
$$
Indeed,
$$
pF^{(k-1)}_{ni}=\sum_{a=-n+1}^{n-1}pF^{(k-2)}_{na}F^{}_{ai}
=\sum_{a=i}^{n-1}pF^{(k-2)}_{na}F^{}_{ai}=
pF^{(k-2)}_{ni}(F^{}_{ii}+n-i-1),
$$
where we have used (2.1) and (5.12). Now (5.13) follows by induction.
Thus, rewriting (5.11) and (5.13) in terms of the generators $z_{ni}$
and $z_{in}$ (see (2.3)) we obtain from (5.10) that
$$
|1-\wt F\ts u^{-1}|_{nn}=1-F^{}_{nn}u^{-1}+
\sum_{i=-n+1}^{n-1}z_{ni}\ts z_{-n,-i}\frac{1}{u(u-f_i-n)}
\prod_{\overset{\ssize a=-n+1}\to {a\ne i}}^{n-1}
\frac{1}{f_i-f_a}.
\tag 5.14
$$

The coefficients of the series $\sdet \Cal F^{(n-1)}(u)$ belong to the
center of $\U(\g_{n-1})$ and so its image under $\pi$ 
coincides with its image
with respect to the Harish-Chandra homomorphism. The latter was 
found in different ways in [Mo3, Section~5] and [MN, Section~6].
The result can be written as follows:
in $\Z(\g_n,\g_{n-1})$ we have
$$
\sdet \Cal F^{(n-1)}(u)=\prod_{a=1}^{n-1}((u-n+3/2)^2-f_a^2)\cdot
\prod_{a=1}^{2n-2} \frac {1}{u-a+1/2}.
\tag 5.15
$$

Finally, using (5.9),
(5.14) and (5.15) we find that in $\Z(\g_n,\g_{n-1})$
$$
s^{\pi}_{n,n}(u)
=c(u)\ts \frac{2u+1}{2u-2n+1}\ts \wh {\Cal F}(-u+n-1)_{nn}=
\frac{u+1/2}{u^{2n}}\ts Z_{n,n}(u).
$$

To find the image of $s_{n,-n}(u)$ under $\pi$ we note that by
(5.1)
$$
[F_{n,-n},s^{\pi}_{n,n}(u)]=-2\ts s^{\pi}_{n,-n}(u).
$$
This implies that the series
$$
\frac{u^{2n}}{u+1/2}\ts s^{\pi}_{n,-n}(u)=
-\frac12\ts [F_{n,-n},Z_{n,n}(u)]
\tag 5.16
$$
is a polynomial in $u$, and by the symmetry relation (3.9)
this polynomial is even. By (5.6) it is of degree $n-1$ in $u^2$
with the highest coefficient $F_{n,-n}$.
Moreover, we see from (5.6) 
that\footnote{If the polynomials 
$Z_{ij}(u)$ are evaluated
in $\R(\h)$, we assume, to avoid
an ambiguity, that they are
written in such a way that the
elements of $\Z(\g_n,\g_{n-1})$ are to the left from 
coefficients belonging to $\R(\h)(u)$, as appears
in their definition.}
$$
Z_{n,n}(-g_i)=-z_{ni}\ts z_{-n,-i},\qquad i=-n+1,\dots,n-1.
$$
Hence, using (2.5) we obtain
$$
[F_{n,-n},Z_{n,n}(-g_i)]=-2\ts z_{ni}\ts z_{n,-i}.
\tag 5.17
$$
Applying the Lagrange interpolation formula to the
polynomial (5.16) by using its values at the $n-1$ points
$-g_i$, $i=1,\dots,n-1$ we prove that
$$
-\frac12\ts [F_{n,-n},Z_{n,n}(u)]=Z_{n,-n}(u).
\tag 5.18
$$

Similarly, replacing $F_{n,-n}$ by $F_{-n,n}$
in the above argument we get
the formula for the image of 
$s_{-n,n}(u)$ under $\pi$.

Finally, to find the image of $s_{-n,-n}(u)$ we use the
following formula implied by (5.1):
$$
[F_{n,-n},[F_{-n,n}, s^{\pi}_{n,n}(u)]]=
2\ts s^{\pi}_{n,n}(u)-2\ts s^{\pi}_{-n,-n}(u). \quad \square
$$

\bigskip
\proclaim
{\bf Theorem 5.2} We have an isomorphism of
$\ts\Y^-(2)$-modules:
$$
V(\lambda)^+_{\mu}\simeq
L(\alpha_1,\beta_1)\ot\cdots\ot L(\alpha_n,\beta_n),
\tag 5.19
$$
where
$$
\aligned
\alpha_i&=\min\{\lambda_{i-1},\mu_{i-1}\}-i+1/2,\qquad i=2,\dots,n,
\qquad \alpha_1=-1/2;\\
\beta_i&=\max\{\lambda_{i},\mu_{i}\}-i+1/2,\qquad i=1,\dots,n-1,
\qquad \beta_n=\lambda_{n}-n+1/2.
\endaligned
\tag 5.20
$$
In particular, $V(\lambda)^+_{\mu}$ is a $\Y(2)$-module.
\endproclaim

\noindent
{\it Proof.}\footnote{Theorem 5.2 was announced in [Mo1].}\ \ 
Let us consider
the following vector $\eta_{\mu}$ in $V(\lambda)^+_{\mu}$: 
$$
\eta_{\mu}=z_{n1'}^{\lambda^{}_{1'}-\mu^{}_{1'}}\cdots
z_{n(n-1)'}^{\lambda^{}_{(n-1)'}-\mu^{}_{(n-1)'}}
\ts\xi, \tag 5.21
$$
where
$$
a'=\cases a,\qquad&\text{if}\quad \lambda_a-\mu_a\geq 0,\\
-a, \qquad&\text{if}\quad \lambda_a-\mu_a< 0,
\endcases
$$
(recall that $\lambda_{-a}=-\lambda_a$, $\ \mu_{-a}=-\mu_a$).
Let us verify that
the vector $\eta_{\mu}$ is nonzero and satisfies the
relations
$$
Z_{-n,n}(u)\ts\eta_{\mu}=0 \tag 5.22
$$
and
$$
Z_{nn}(u)\ts\eta_{\mu}=(u-\alpha_2)\cdots 
(u-\alpha_n)(u+\beta_1)\cdots (u+\beta_n)\ts\eta_{\mu}.
\tag 5.23
$$
Using (2.6) and (2.7) we easily obtain by
induction on the degree of the monomial (5.21) that
$$
z_{-i',n}\ts \eta_{\mu}=0,\qquad i=1,\dots,n-1.
\tag 5.24
$$
Let us verify that for $i=1,\dots,n-1$
$$
z_{in}\ts\eta_{\mu}=-(m_i+\alpha_2)\cdots 
(m_i+\alpha_n)(m_i-\beta_1)\cdots (m_i-\beta_n)\ts \eta_{\mu+\delta_i},
\tag 5.25
$$
if $\lambda_i\geq\mu_i$; and
$$
z_{-i,n}\ts\eta_{\mu}=-(m_i-\alpha_2-1)\cdots 
(m_i-\alpha_n-1)(m_i+\beta_1-1)\cdots (m_i+\beta_n-1)
\ts \eta_{\mu-\delta_i},
\tag 5.26
$$
if $\lambda_i\leq\mu_i$; here we have set 
$m_i=\mu_i-i+1/2$. We shall prove (5.23), (5.25) and (5.26)
simultaneously by induction on the degree of the monomial (5.21).
If the degree is zero all the relations are obvious.
If $\lambda_i=\mu_i$ then both (5.25) and (5.26) hold
by (5.24) because
in this case $\beta_i=\alpha_{i+1}+1=m_i$. 
Suppose now that $\lambda_i>\mu_i$.
By (2.6) we can write
$$
z_{in}\ts\eta_{\mu}=z_{in}\ts z_{ni}\ts \eta_{\mu+\delta_i}.
$$
Formula (5.7) gives $z_{in}\ts z_{ni}=-Z_{-n,-n}(-g_{-i})$. However,
$$
(-g_{-i})\ts \eta_{\mu+\delta_i}=m_i\ts \eta_{\mu+\delta_i}.
$$
Applying Theorem 5.1 and the symmetry relation (3.9) we
obtain
$$
Z_{-n,-n}(m_i)=\left(\frac{1}{2m_i}-1\right)Z_{n,n}(-m_i)
-\frac{1}{2m_i}\ts Z_{n,n}(m_i).
$$
Using the induction hypotheses we can find 
$Z_{n,n}(\pm m_i)\ts \eta_{\mu+\delta_i}$ by (5.23).
Since $m_i=\alpha_{i+1}$ we have
$Z_{n,n}(m_i)\ts \eta_{\mu+\delta_i}=0$ and (5.25)
follows.

The same argument can be applied to prove (5.26). Here we use
the relation $z_{-i,n}\ts z_{n,-i}=Z_{-n,-n}(-g_i)$
implied by (5.7).

To prove (5.23) it suffices to check this relation for $2n-1$ different
values of $u$ because $Z_{n,n}(u)$ is a monic
polynomial in $u$ of degree $2n-1$. For these take
the eigenvalues of 
$-g_i$ with $i=-n+1,\dots,n$ 
on the vector $\eta_{\mu}$
and then use (5.6), (5.25) and (5.26).

Relation (5.22) follows from (2.5) and (5.24)--(5.26). The fact that
$\eta_{\mu}\ne 0$ is implied by (5.25) and (5.26). Indeed,
applying appropriate operators $z_{in}$ or $z_{-i,n}$
to the vector $\eta_{\mu}$
repeatedly we can obtain the highest weight vector
$\xi$ of $V(\lambda)$ with a nonzero coefficient.

Finally, using formulas (5.8) we deduce from
(5.22) and (5.23) that $\eta_{\mu}$
is the highest vector of the representation $V(\lambda)^+_{\mu}$
of the algebra $\Y^-(2)$, and the highest weight
is given by
$$
\mu(u)=(1-\alpha_1u^{-1})\cdots (1-\alpha_nu^{-1})
(1+\beta_1u^{-1})\cdots (1+\beta_nu^{-1}).
$$
On the other hand, Proposition 4.1 implies that the representation
(5.19) of $\Y^-(2)$ is irreducible because
the strings
$$
S(\alpha_1,\beta_1),\dots,S(\alpha_n,\beta_n),S(-\beta_1,-\alpha_1),
\dots, S(-\beta_n,-\alpha_n)
$$
are pairwise in general position.
Therefore, it is isomorphic
to the highest weight representation $V(\mu(u))\simeq V(\lambda)^+_{\mu}$;
see (4.5). $\square$
\medskip

The branching rule for representations of the symplectic
Lie algebras [Z1] follows immediately from Theorem~5.2.

\bigskip
\proclaim
{\bf Corollary 5.3} 
The restriction of $V(\lambda)$ to the subalgebra $\g_{n-1}$
is isomorphic to the direct sum (1.3) of irreducible finite-dimensional
representations $V'(\mu)$ of $\g_{n-1}$
where the multiplicity $c(\mu)$ equals the number of $n$-tuples of
integers $\nu$ satisfying the inequalities (1.4).
\endproclaim

\Proof Examining the weight diagram of $V(\lambda)$
(see, e.g., [D, Section 7.2]) we find that
the representation
$V'(\mu)$ can occur in the decomposition (1.3) only if
the $\g_{n-1}$-highest weight $\mu$ satisfies the condition
$$
\lambda_{i+1}\leq\mu_i\leq\lambda_{i-1},\qquad i=1,\dots,n-1,
$$
($\lambda_0:=0$). In this case the multiplicity $c(\mu)$ 
coincides with the dimension
of the space $V(\lambda)^+_{\mu}$. By Theorem~5.2,
$$
\dim V(\lambda)^+_{\mu}=\prod_{i=1}^n(\alpha_i-\beta_i+1),
$$
which is equal to the number of solutions of the inequalities (1.4).
$\square$

\bigskip
\bigskip
\noindent
{\bf 6. Construction of the basis in $V(\lambda)$}
\bigskip

In this section we complete the proof of Theorem 1.1.

We first construct a basis in the space $V(\lambda)^+_{\mu}$.
A basis in $V(\lambda)$ will then be obtained by induction
with the use of the branching rule (1.3).

Note that by (2.12), (5.6) and (5.18) we
have $Z_{n,-n}(g_n)=z_{n,-n}$ and so
the polynomial $Z_{n,-n}(u)$ can be written as follows:
$$
Z_{n,-n}(u)=
\sum_{i=1}^{n}z_{ni}\ts z_{n,-i}\prod_{a=1,\ts a\ne i}^{n}
\frac{u^2-g_a^2}{g_i^2-g_a^2},
$$
where, as before, $z_{nn}=1$ and $g_i=F_{ii}-i+1/2$.

We shall use the notation
$$
l_i=\lambda_i-i+1/2,\qquad \gamma_i=\nu_i-i+1/2,
\qquad m_i=\mu_i-i+1/2
\tag 6.1
$$
with $i$ ranging over $\{1,\dots,n\}$ or $\{1,\dots,n-1\}$ respectively.
Given $\nu$ satisfying the inequalities (1.4) consider the vector
$$
\xi_{\nu}=\prod_{i=1}^{n-1}
z_{ni}^{\nu_i-\mu_i}\ts z_{n,-i}^{\nu_i-\lambda_i}\cdot
\prod_{k=l_n}^{\gamma_n-1}
Z_{n,-n}(k)\ts\xi,
\tag 6.2
$$
where the operators $z_{nj}$ and $Z_{n,-n}(u)$ are defined in
(2.3) and (5.4). 
Here the action of elements of $\R(\h)$
is obtained by the extension from that of $\U(\h)$. 
The action is well-defined
for those elements whose denominators are not zero operators.
One easily checks that the vector $\xi_{\nu}$ is well-defined.
By definition of the algebra $\Z(\g_n,\g_{n-1})$ we have
$\xi_{\nu}\in V(\lambda)^+$. Moreover, $\xi_{\nu}$
is clearly of $\g_{n-1}$-weight $\mu$.

\bigskip
\proclaim
{\bf Proposition 6.1} The vectors $\xi_{\nu}$ 
with $\nu$ satisfying the inequalities
(1.4) form a basis in the space $V(\lambda)^+_{\mu}$.
\endproclaim

\Proof We use Theorem 5.2. The strings $S(\alpha_i,\beta_i)$
obviously satisfy the condition (4.10).
By Proposition 4.2 the vectors $\eta_{\gamma}$ defined by (4.9)
with the $n$-tuple $\gamma=(\gamma_i)$ given in (6.1)
form a basis in $V(\lambda)^+_{\mu}$. However, the
operator $S^{\natural}_{n,-n}(u)$ coincides with $Z_{n,-n}(u)$
by Theorem 5.1.
That is, the vectors
$$
\eta_{\gamma}=\prod_{i=1}^n Z_{n,-n}(\gamma_i-1)\cdots 
Z_{n,-n}(\beta_i+1)Z_{n,-n}(\beta_i)\ts\eta_{\mu},
\tag 6.3
$$
with $\eta_{\mu}$ defined in (5.21) form a basis in
$V(\lambda)^+_{\mu}$. Let us show that
for each $\nu$ satisfying (1.4) we have the equality of
corresponding vectors:
$$
\eta_{\gamma}=\xi_{\nu}.\tag 6.4
$$
Note first that
for any $i=-n+1,\dots,n-1$ and any value
$u\in\C$ one has
$$
[Z_{n,-n}(u), z_{ni}]=0. \tag 6.5
$$
Indeed, by (5.1) and (5.8) we have $[Z_{n,-n}(u),F_{ni}]=0$.
It remains to apply the extremal projection $p=p_{n-1}$ 
and use the fact that $Z_{n,-n}(u)$ commutes with $\g_{n-1}$.

Let $b_1,\dots,b_k$ be all the indices $a$ among
$1,\dots,n-1$ for which the difference $\lambda_a-\mu_a$
is positive and let $c_1,\dots,c_l$ be the remaining indices; $k+l=n-1$.
Using (2.6) and (6.5) rewrite the vector 
(6.3) as follows:
$$
\eta_{\gamma}=\prod_{i=1}^kz_{nb_i}^{\lambda_{b_i}-\mu_{b_i}}
{}\cdot\prod_{i=1}^n Z_{n,-n}(\gamma_i-1)\cdots 
Z_{n,-n}(\beta_i)\cdot
\prod_{i=1}^lz_{n,-c_i}^{\mu_{c_i}-\lambda_{c_i}}\ts\xi.
\tag 6.6
$$
Further, by (5.4) we have
$$
Z_{n,-n}(g_i)=z_{ni}\ts z_{n,-i},\qquad i=1,\dots, n-1.
$$
However,
$$
g_i\ts\prod_{i=1}^lz_{n,-c_i}^{\mu_{c_i}-\lambda_{c_i}}\ts\xi
=\beta_i\ts \prod_{i=1}^lz_{n,-c_i}^{\mu_{c_i}-\lambda_{c_i}}\ts\xi.
$$
Therefore, given $i\in\{1,\dots,n-1\}$ the operator
$Z_{n,-n}(\beta_i)$ in (6.6) 
can be replaced with $z_{ni}\ts z_{n,-i}$.
Moving the $z_{ni}$ to the left permuting it with the operators
of the form $Z_{n,-n}(u)$ we represent the vector again in the form
(6.6). Proceeding by induction we shall get the expression for the vector
$\eta_{\gamma}$ which coincides with (6.2).
$\square$
\medskip

Given a pattern $\Lambda$ associated with $\lambda$ (see Section 1)
define the vector
$\xi^{}_{\Lambda}$ by the formula
$$
\xi^{}_{\Lambda}=\prod_{k=1,\dots,n}^{\rightarrow}
\left(\prod_{i=1}^{k-1}z_{ki}^{\lambda'_{ki}-\lambda^{}_{k-1,i}}
z_{k,-i}^{\lambda'_{ki}-\lambda^{}_{ki}}\cdot
\prod_{a=l^{}_{kk}}^{l'_{kk}-1}
Z_{k,-k}(a+1/2)\right)\xi.
$$
Here the $z_{kj}$ and $Z_{k,-k}(u)$ are elements of the Mickelsson
algebra $\Z(\g_k,\g_{k-1})$,
and the numbers $l^{}_{kk}, l'_{kk}$ are defined in (1.5). 
The branching rule (1.3) and Proposition~6.1
immediately imply the following.

\bigskip
\proclaim
{\bf Proposition 6.2} The vectors $\xi^{}_{\Lambda}$ where $\Lambda$
runs over all patterns associated with $\lambda$ form a basis
in the representation $V(\lambda)$ of $\g_n$. $\square$
\endproclaim

\medskip

Our next task is to calculate the matrix elements of the generators
$F_{kk}$, $F_{k,-k}$, $F_{-k,k}$, $F_{k-1,-k}$ 
of the Lie algebra $\g_n$ in the basis 
$\{\xi^{}_{\Lambda}\}$. Note that the elements 
$F_{kk}$, $F_{k,-k}$, $F_{-k,k}$ belong to the centralizer 
of the subalgebra $\g_{k-1}$ in $\U(\g_k)$. Therefore, these
operators preserve the subspace of $\g_{k-1}$-highest vectors
in $V(\lambda)$. So, it suffices to compute the action
of these operators with $k=n$ in the basis $\{\xi_{\nu}\}$
of the space $V(\lambda)^+_{\mu}$; see Proposition 6.1.

For $F_{nn}$ we immediately get
$$
F_{nn}\ts\xi_{\nu}=\left(2\ts\sum_{i=1}^n\nu_i-
\sum_{i=1}^n\lambda_i-\sum_{i=1}^{n-1}\mu_i\right)\xi_{\nu}.
\tag 6.7
$$

Further, by (6.3) and (6.4)
$$
Z_{n,-n}(\gamma_i)\ts\xi_{\nu}=\xi_{\nu+\delta_i},\qquad i=1,\dots,n.
$$
However, $Z_{n,-n}(u)$ is a polynomial in $u^2$ of degree $n-1$
with the highest coefficient $F_{n,-n}$; see (5.4). Applying 
the Lagrange interpolation
formula with the interpolation points 
$\gamma_i$, $i=1,\dots,n$
we obtain
$$
Z_{n,-n}(u)\ts\xi_{\nu}=\sum_{i=1}^n\prod_{a=1,\ts a\ne i}^n
\frac{u^2-\gamma_a^2}{\gamma_i^2-\gamma_a^2}\ \xi_{\nu+\delta_i}.
$$
Taking here the coefficient at $u^{2n-2}$ we get
$$
F_{n,-n}\ts\xi_{\nu}=\sum_{i=1}^n\prod_{a=1,\ts a\ne i}^n
\frac{1}{\gamma_i^2-\gamma_a^2}\ \xi_{\nu+\delta_i}.
\tag 6.8
$$

Similarly, we see from (5.5) that 
$Z_{-n,n}(u)$ is a polynomial in $u^2$ of degree $n-1$
with the highest coefficient $F_{-n,n}$. Using the defining relations
(3.1) we can write the following formula for the operator
$S_{-n,n}(u)$ in $V(\lambda)^+_{\mu}$:
$$
S_{-n,n}(u)=\frac{(-1)^n(u+1/2)}{u}\left(T_{-n,n}(u)T_{-n,-n}(-u)-
T_{-n,n}(-u)T_{-n,-n}(u)\right).
$$
Hence, by Theorem 5.1 we obtain the equality
of operators in $V(\lambda)^+_{\mu}$:
$$
Z_{-n,n}(u)=\frac{(-1)^n}{u}\left(T_{-n,n}(u)T_{-n,-n}(-u)-
T_{-n,n}(-u)T_{-n,-n}(u)\right).
$$
This implies that $F_{-n,n}$,
as an operator in $V(\lambda)^+_{\mu}$,
coincides with $2\ts t_{-n,n}^{(1)}$, where $t_{-n,n}^{(1)}$
is the highest coefficient
of the polynomial $T_{-n,n}(u)$ which has degree $n-1$; see (4.6).
Using (4.13) we find that
$$
T_{-n,n}(-\gamma_i)\ts \xi_{\nu}=
-2\ts\prod_{k=1}^n(\alpha_k-\gamma_i+1)(\beta_k-\gamma_i)\cdot
\prod_{a=1,\ts a\ne i}^n (-\gamma_i-\gamma_a+1)\ts \xi_{\nu-\delta_i}.
$$
Note that by (5.20) we have
$$
\prod_{k=1}^n(\alpha_k-\gamma_i+1)(\beta_k-\gamma_i)
=(1/2-\gamma_i)\prod_{k=1}^n(l_k-\gamma_i)\prod_{k=1}^{n-1}(m_k-\gamma_i);
$$
see (6.1).
Applying the Lagrange interpolation formula to the polynomial
$T_{-n,n}(u)$ with the interpolation points $-\gamma_i$,
$i=1,\dots,n$
and taking the coefficient at $u^{n-1}$
we finally obtain that
$$
F_{-n,n}\ts \xi_{\nu}=
2\ts\sum_{i=1}^n\ 
\frac{\dsize\prod_{a=1}^n(l_a-\gamma_i)(\gamma_a+\gamma_i-1)
\prod_{a=1}^{n-1}(m_a-\gamma_i)}
{\dsize\prod_{a=1,\ts a\ne i}^n 
(\gamma_i-\gamma_a)}\  \xi_{\nu-\delta_i}.
\tag 6.9
$$
\medskip

To compute the action of the elements $F_{k-1,-k}$ we may only
consider the case $k=n$. The operator $F_{n-1,-n}$ preserves
the subspace of $\g_{n-2}$-highest vectors in $V(\lambda)$.
Therefore it suffices to calculate its action
on the basis vectors of the form
$$
\xi_{\nu\mu\nu'}=\prod_{i=1}^{n-2}
z_{n-1,i}^{\nu'_i-\mu'_i}\ts z_{n-1,-i}^{\nu'_i-\mu_i}\cdot
\prod_{a=m_{n-1}}^{\gamma'_{n-1}-1}
Z_{n-1,-n+1}(a)\ts\xi_{\nu\mu},
$$
where $\xi_{\nu\mu}=\xi_{\nu}$ is defined in (6.2),
$\mu'$ is a fixed 
$\g_{n-2}$-highest weight, $\nu'$ is an $(n-1)$-tuple of integers
such that the inequalities (1.4)
are satisfied with $\lambda$, $\nu$, $\mu$ respectively replaced
by $\mu$, $\nu'$, $\mu'$, and we set $\gamma'_i=\nu'_i-i+1/2$. 
The operator $F_{n-1,-n}$ is permutable with the 
elements $z_{n-1,i}$ and
$Z_{n-1,-n+1}(u)$ which follows from the explicit formulas
(2.3) and (5.4). Hence, we can write
$$
F_{n-1,-n}\ts \xi_{\nu\mu\nu'}=X_{\mu\nu'}\ts F_{n-1,-n}\ts \xi_{\nu\mu},
$$
where $X_{\mu\nu'}$ denotes the operator
$$
X_{\mu\nu'}=\prod_{i=1}^{n-2}
z_{n-1,i}^{\nu'_i-\mu'_i}\ts z_{n-1,-i}^{\nu'_i-\mu_i}\cdot
\prod_{a=m_{n-1}}^{\gamma'_{n-1}-1}
Z_{n-1,-n+1}(a).
$$
Let us now apply formula (2.13). We have
$$
f_a\ts \xi_{\nu\mu}=(m_a-1/2)\ts \xi_{\nu\mu},\qquad a=1,\dots,n-1;
\tag 6.10
$$
see (6.2). Further, for $i=1,\dots,n-1$
$$
X_{\mu\nu'}\ts z_{n-1,-i}z_{ni}\ts
\xi_{\nu\mu}=\xi_{\nu,\mu-\delta_i,\nu'}.
\tag 6.11
$$
Indeed, by (6.2)
$
z_{ni}\ts \xi_{\nu\mu}=\xi_{\nu,\mu-\delta_i},
$
and
$
X_{\mu\nu'}\ts z_{n-1,-i}=X_{\mu-\delta_i,\nu'}
$
for $i<n-1$,
where we have used (2.6) and (6.5). Let us verify that
the latter formula holds for $i=n-1$ as well.
By (2.12) and (5.4)
we can write
$
z_{n-1,-n+1}=Z_{n-1,-n+1}(g_{n-1}).
$
However, 
$$
g_{n-1}\ts
\xi_{\nu,\mu-\delta_{n-1}}=(m_{n-1}-1)\ts 
\xi_{\nu,\mu-\delta_{n-1}},
$$
and so
$$
X_{\mu\nu'}\ts z_{n-1,-n+1}=X_{\mu\nu'}Z_{n-1,-n+1}(m_{n-1}-1)
=X_{\mu-\delta_{n-1},\nu'},
$$
as desired.

Finally, for $j=1,\dots,n-1$ consider the expression
$$
X_{\mu\nu'}\ts z_{n-1,j}\ts z_{n,-j}\ts \xi_{\nu\mu}.
\tag 6.12
$$
First transform the vector $z_{n,-j}\ts \xi_{\nu}$. The calculation
is trivial if $\nu_j=\mu_j$, so
we shall assume that $\nu_j-\mu_j\geq 1$.
We have
$$
z_{n,-j}\ts \xi_{\nu\mu}=z_{n,-j}\ts z_{nj}\ts\xi_{\nu,\mu+\delta_j}.
$$
By (5.17) and (5.18) we have $z_{n,-j}\ts
z_{nj}=Z_{n,-n}(g_j-1)$. However,
$$
(g_j-1)\ts \xi_{\nu,\mu+\delta_j}=m_j\ts 
\xi_{\nu,\mu+\delta_j}.
$$
To calculate $Z_{n,-n}(m_j)\xi_{\nu,\mu+\delta_j}$
we apply again the Lagrange interpolation formula (cf. the proof
of (6.8))
for the polynomial $Z_{n,-n}(u)$ at the interpolation points
$\gamma_i$, $i=1,\dots,n$ and then put $u=m_j$. The result is
$$
Z_{n,-n}(m_j)\xi_{\nu,\mu+\delta_j}=
\sum_{i=1}^n \prod_{a=1,\ts a\ne i}^n\frac{m_j^2-\gamma_a^2}
{\gamma_i^2-\gamma_a^2}\ 
\xi_{\nu+\delta_i,\mu+\delta_j}.
\tag 6.13
$$

Let us now transform the operator $X_{\mu\nu'}\ts z_{n-1,j}$, 
$j<n-1$.
Here the calculation is very similar to the previous one.
We shall assume that $\nu'_j-\mu_j\geq 1$. We have
$$
X_{\mu\nu'}\ts z_{n-1,j}=X_{\mu+\delta_j,\nu'}\ts z_{n-1,-j}\ts z_{n-1,j}.
$$
Using (5.17) and (5.18) again, we can write 
$z_{n-1,-j}\ts z_{n-1,j}=Z_{n-1,-n+1}(g_j-1)$. We have
$$
(g_j-1)\ts\xi_{\nu+\delta_i,\mu+\delta_j}=m_j
\ts\xi_{\nu+\delta_i,\mu+\delta_j}.
$$
Exactly as above, we use the Lagrange interpolation formula
for the polynomial $Z_{n-1,-n+1}(u)$ with the interpolation points
$\gamma'_r$, $r=1,\dots,n-1$ and then put $u=m_j$. This gives
$$
X_{\mu+\delta_j,\nu'}Z_{n-1,-n+1}(m_j)=
\sum_{r=1}^{n-1} \prod_{a=1,\ts a\ne r}^{n-1}\frac{m_j^2-{\gamma'_a}^2}
{{\gamma'_r}^2-{\gamma'_a}^2}\ts
X_{\mu+\delta_j,\nu'+\delta_r}.
\tag 6.14
$$
In the case $j=n-1$ we write
$X_{\mu\nu'}=X_{\mu+\delta_{n-1},\nu'}Z_{n-1,-n+1}(m_{n-1})$
and (6.14) holds for this case as well.

Combining (6.10)--(6.14) we obtain from (2.13)
$$
\align
F_{n-1,-n}\ts \xi_{\nu\mu\nu'}=\sum_{i=1}^{n-1}\frac{1}{2m_i-1}
&\prod_{a=1,\ts a\ne i}^{n-1}\frac {1}{(m_i-m_a)(m_i+m_a-1)}\ 
\xi_{\nu,\mu-\delta_i,\nu'}\\
{}+\sum_{i=1}^{n}\sum_{j,r=1}^{n-1}\frac{1}{2m_j-1}
&\prod_{a=1,\ts a\ne j}^{n-1}\frac {1}{(m_j-m_a)(m_j+m_a-1)}\tag 6.15\\
&\prod_{a=1,\ts a\ne i}^n\frac{m_j^2-\gamma_a^2}
{\gamma_i^2-\gamma_a^2}
\prod_{a=1,\ts a\ne r}^{n-1}\frac{m_j^2-{\gamma'_a}^2}
{{\gamma'_r}^2-{\gamma'_a}^2}\ 
\xi_{\nu+\delta_i,\mu+\delta_j,\nu'+\delta_r}.
\endalign
$$

To complete the proof of Theorem 1.1 we rewrite the formulas
(6.7)--(6.9) and (6.15) in the notation related to the
patterns $\Lambda$ (see Section 1) to get the matrix 
elements of the generators $F_{kk}$, $F_{k,-k}$, $F_{-k,k}$, $F_{k-1,-k}$
in the basis $\{\xi^{}_{\Lambda}\}$ of $V(\lambda)$. 
The formulas of Theorem~1.1 are given in the normalized basis
$\{\zeta^{}_{\Lambda}\}$
where
$$
\zeta^{}_{\Lambda}=N^{}_{\Lambda}\ts \xi^{}_{\Lambda},
$$
and
$$
N^{}_{\Lambda}=\prod_{k=2}^n\prod_{1\leq i<j\leq k}(-l'_{ki}-l'_{kj}-1)!
$$

Theorem 1.1 is proved.

\bigskip
\noindent
{\it Remark 6.3.\ }  As it follows from our arguments, 
the problem of constructing
an orthogonal basis in the $\g_n$-module
$V(\lambda)$ is in fact reduced to the problem of constructing
an orthogonal basis in the $\Y^-(2)$-module
$V(\lambda)^+_{\mu}$. A natural way to do this
is to find the eigenvectors of a commutative subalgebra
in $\Y^-(2)$. Such subalgebra is
generated by the coefficients of the series 
$s_{-1,-1}(u)+s_{1,1}(u)=\tr S(u)$.
The problem can also be
reformulated for the Yangian $\Y(2)$ since the commutative subalgebra
can be identified with that of $\Y(2)$.

\bigskip
\bigskip
\noindent
{\bf References}
\bigskip

\itemitem{[AST1]}
{Asherova, R. M., Smirnov, Yu. F., Tolstoy, V. N.:}
{Projection operators for simple Lie groups}.
{Theor. Math. Phys.} {\bf 8}, 813--825 (1971)

\itemitem{[AST2]}
{Asherova, R. M., Smirnov, Yu. F., Tolstoy, V. N.:}
{Projection operators for simple Lie groups. II. General scheme for
constructing lowering operators. The groups 
$SU(n)$}. 
{Theor. Math. Phys.} {\bf 15}, 392--401 (1973)

\itemitem{[AST3]}
{Asherova, R. M., Smirnov, Yu. F., Tolstoy, V. N.:}
{Description 
of a certain class of projection operators for complex
semisimple Lie algebras}. 
Math. Notes {\bf 26}, no. 1-2, 499 -- 504  (1979)

\itemitem{[B]}
{Berele, A.:} {Construction of $\text{ Sp}$-modules by
   tableaux}. Linear and Multilinear Algebra 
{\bf 19}, 299--307  (1986)

\itemitem{[Bi1]}
{Bincer, A.:} {Missing label operators in the reduction
$Sp(2n)\downarrow Sp(2n-2)$}.
J. Math. Phys. {\bf 21}, 671--674  (1980)

\itemitem{[Bi2]}
{Bincer, A.:}  {Mickelsson lowering operators 
for the symplectic group}.
J. Math. Phys. {\bf 23}, 347--349  (1982)

\itemitem{[CP]}
{Chari, V., Pressley, A.:}
{Yangians and $R$-matrices}.
{L'Enseign. Math.}
{\bf 36},
267--302 (1990)

\itemitem{[C]}
{Cherednik, I. V.:}
{A new interpretation of Gelfand--Tzetlin bases}. {Duke Math. J.}
{\bf 54},
563--577 (1987)

\itemitem{[CK]}
{De Concini, C., Kazhdan, D.:} {Special bases
   for $S_{N}$ and $\text{GL}(n)$}. 
Israel J. Math. {\bf 40}, no. 3-4,
   275--290  (1981)

\itemitem{[D]}
{Dixmier, J.:}
{Alg\`ebres Enveloppantes}. 
{Paris: Gauthier-Villars}
1974

\itemitem{[Dr]}
{Drinfeld, V. G.:}
{A new realization of Yangians and quantized affine
algebras}.
{Soviet Math. Dokl.}
{\bf 36},
212--216 (1988)

\itemitem{[GKLLRT]} 
{Gelfand, I. M., Krob, D., Lascoux, A.,
Leclerc, B., Retakh, V. S., Thibon, J.-Y.:}
{Noncommutative symmetric functions}. Adv. Math. 
{\bf 112}, 218--348 (1995)

\itemitem{[GT1]}
{Gelfand, I. M., Tsetlin, M. L.:}
{Finite-dimensional representations of the group
of unimodular matrices}.
Dokl. Akad. Nauk SSSR {\bf 71}, 825--828  (1950) (Russian).
English transl. in: Gelfand, I. M. Collected papers. Vol II,
Berlin: Springer-Verlag  1988

\itemitem{[GT2]}
{Gelfand, I. M., Tsetlin, M. L.:}
{Finite-dimensional representations of groups of
orthogonal matrices}.
Dokl. Akad. Nauk SSSR {\bf 71}, 1017--1020  (1950) (Russian).
English transl. in: Gelfand, I. M. Collected papers. Vol II,
Berlin: Springer-Verlag  1988

\itemitem{[GZ1]}
{Gelfand, I. M., Zelevinsky, A.:} 
{Models of representations of classical
groups and their hidden symmetries}. Funct. Anal. Appl. {\bf 18},
183--198  (1984)

\itemitem{[GZ2]}
{Gelfand, I. M., Zelevinsky, A.:}  
{Multiplicities and
   proper bases for $\text{gl}_n$}. In: Group theoretical methods in
   physics. Vol. II, Yurmala 1985, pp. 147--159.
Utrecht: VNU Sci. Press 1986

\itemitem{[G1]}
{Gould, M. D.:}
{On the matrix elements of 
the $\text{U}(n)$ generators}. 
J. Math. Phys. {\bf 22}, 15--22  (1981)

\itemitem{[G2]}
{Gould, M. D.:} 
{Wigner coefficients for a semisimple
Lie group and the matrix elements of the $\text {O}(n)$ generators}. J.
Math. Phys. {\bf 22}, 2376--2388  (1981)

\itemitem{[G3]}
{Gould, M. D.:} 
 {Representation theory of the symplectic
groups. I}. J. Math. Phys. {\bf 30}, 1205--1218  (1989)

\itemitem{[GK]}
{Gould, M. D., Kalnins, E. G.:}  
{A projection-based
   solution to the $\text {Sp}(2n)$ state labeling problem}. J. Math. Phys.
   {\bf 26}, 1446--1457  (1985)

\itemitem{[H]}
{Hegerfeldt, G. C.:}
{Branching theorem for the symplectic groups}.
{J. Math. Phys.} {\bf 8}, 1195--1196  (1967)

\itemitem{[Ho]}
{Van den Hombergh, A.:}
{A note on Mickelsson's step algebra}.
{Indag. Math.} {\bf 37}, no.1, 42--47  (1975)

\itemitem{[Ka]}
{Kashiwara, M.:}
 {Crystalizing the $q$-analogue of
   universal enveloping algebras}. Comm. Math. Phys. {\bf 133},
   249--260  (1990)

\itemitem{[Ki]}
{King, R. C.:}
{Weight multiplicities for the classical
   groups}. In: Group theoretical methods in physics. Fourth Internat.
Colloq., Nijmegen 1975. Lecture Notes in Phys., Vol.
   50, pp. 490--499. Berlin: Springer 1976

\itemitem{[KS]}
{King, R. C., El-Sharkaway, N. G. I.:}
{Standard
   Young tableaux and weight multiplicities of the classical Lie groups}.
   J. Phys. A {\bf 16}, 3153--3177  (1983)

\itemitem{[KW]}
{King, R. C., Welsh, T. A.:} {Construction of orthogonal
group modules using tableaux}. Linear and Multilinear Algebra 
{\bf 33}, 251--283  (1993)

\itemitem{[K]}
{Kirillov, A. A.:}
{A remark on the Gelfand-Tsetlin patterns for symplectic groups}.
J. Geom. Phys.
   {\bf 5}, 473--482  (1988)

\itemitem{[KT]}
{Koike, K., Terada, I.:}
 {Young-diagrammatic
   methods for the representation theory of the classical groups of type
   $B_n,\ts C_n,\ts D_n$}. J. Algebra {\bf 107}, 466--511
 (1987)

\itemitem{[LMS]}
{Lakshmibai, V., Musili, C., Seshadri, C. S.:}
{Geometry of $G/P$. IV. Standard monomial theory for classical types}.
Proc. Indian Acad. Sci. Sect. A Math. Sci. {\bf 88}, no. 4,
279--362  (1979)

\itemitem{[L]}
{Littelmann, P.:}
{An algorithm to compute bases and
   representation matrices for ${SL}_{n+1}$-representations}.
J. Pure Appl. Algebra
   {\bf 117}/{\bf 118}, 447--468  (1997)

\itemitem{[Lu]}
{Lusztig, G.:}
 {Canonical bases arising from quantized
   enveloping algebras}. J. Amer. Math. Soc. {\bf 3}, 
447--498  (1990)

\itemitem{[M1]}
{Mathieu, O.:}
 {Good bases for $G$-modules}. Geom.
Dedicata {\bf 36}, 51--66  (1990)

\itemitem{[M2]}
{Mathieu, O.:}
{Bases\ \  des\ \  repr\'esentations\ \  des\ \ 
   groupes simples complexes (d'apr\`es Kashiwara, Lusztig, Ringel et
al.)}.
S\'emin. Bourbaki, Vol. 1990/91.
   Ast\'erisque no. 201--203. Exp. no. 743, 421--442  (1992)

\itemitem{[Mi1]} 
{Mickelsson, J.:}
{Lowering operators and the
symplectic group}. Rep. Math. Phys. {\bf 3}, 193--199  (1972)

\itemitem{[Mi2]}
{Mickelsson, J.:}
{Step algebras of semi-simple
subalgebras of Lie algebras}. Rep. Math. Phys. {\bf 4},
307--318  (1973)

\itemitem{[Mo1]}
{Molev, A.:}
{Representations of twisted Yangians}. {Lett. Math. Phys.}
{\bf 26},
211--218 (1992)

\itemitem{[Mo2]}
{Molev, A.:}
{Gelfand--Tsetlin basis for representations of Yangians}.
{Lett. Math. Phys.}
{\bf 30},
53--60 (1994)

\itemitem{[Mo3]}
{Molev, A.:}
{Sklyanin determinant, Laplace operators, and characteristic identities
for classical Lie algebras}. 
{J. Math. Phys.}
{\bf 36},
923--943 (1995)

\itemitem{[Mo4]} 
{Molev, A.:}
{Finite-dimensional irreducible representations of twisted
Yangians}.
Preprint CMA 047-97, Austral. Nat. University, Canberra;
q-alg/9711022.

\itemitem{[MN]} 
{Molev, A., Nazarov, M.:}
{Capelli identities for classical Lie algebras}.
Preprint CMA 003-97, Austral. Nat. University, Canberra;
q-alg/9712021.

\itemitem{[MNO]}
{Molev, A., Nazarov, M., Olshanski, G.:}
{ Yangians and classical Lie algebras}.
Russian Math. Surveys
{\bf 51}:2,
205--282 (1996)

\itemitem{[MO]}
{Molev, A., Olshanski, G.:}
{Centralizer construction for twisted Yangians}.
Preprint CMA 065-97, Austral. Nat. University, Canberra;
q-alg/9712050.

\itemitem{[NM]}
{Nagel, J. G., Moshinsky, M.:}
 {Operators
that lower or raise the irreducible vector spaces of $U_{n-1}$
contained in an irreducible vector space of $U_n$}.
{J. Math. Phys.} {\bf 6}, 682--694  (1965)

\itemitem{[NT1]}
{Nazarov, M., Tarasov, V.:}
{Yangians and Gelfand--Zetlin bases}. Publ. RIMS, Kyoto Univ.
{\bf 30}, 459--478  (1994)

\itemitem{[NT2]}
{Nazarov, M., Tarasov, V.:}
{Representations of Yangians with Gelfand--Zetlin bases}. 
{J. Reine Angew. Math.} {\bf 496}, 181--212  (1998)

\itemitem{[Ok]}
{Okounkov, A.:}
{Multiplicities and Newton polytopes}.
In: Olshanski,~G. (ed.)
Kirillov's Seminar on Representation Theory.
{Amer. Math. Soc. Transl.}
{\bf 181},
 pp. 231--244. AMS,
Providence RI 1998

\itemitem{[O1]} 
{Olshanski, G. I.:}
{Extension of the algebra $U(g)$ for 
infinite-dimensional classical Lie algebras $g$, 
and the Yangians $Y(gl(m))$.}
{Soviet Math. Dokl.} {\bf 36}, 569--573  (1988)

\itemitem{[O2]}
{Olshanski, G. I.:}
{Representations of 
infinite-dimensional classical groups, limits of enveloping algebras, and
Yangians}. In: Kirillov,~A.~A. (ed.)
Topics in Representation Theory.
{Advances in Soviet Math.} {\bf 2}, pp.
1--66. AMS, Providence RI  1991

\itemitem{[O3]}
{Olshanski, G. I.:}
{Twisted Yangians and infinite-dimensional classical Lie algebras}.
In: Kulish,~P.~P. (ed.)
Quantum Groups, {Lecture Notes in Math.}
{\bf 1510}, pp. 103--120.
Berlin-Heidelberg: Springer
1992

\itemitem{[PH]}
{Pang, S. C., Hecht, K. T.:}
{Lowering and raising
operators for the orthogonal group in the chain ${O}(n)\supset
{O}(n-1)\supset \cdots $, and their graphs}.
{J. Math. Phys.}
{\bf 8}, 1233--1251  (1967)

\itemitem{[P1]}
{Proctor, R.:}
 {Odd symplectic groups}.
Invent.
   Math. {\bf 92}, 307--332  (1988)

\itemitem{[P2]}
{Proctor, R.:}
{Young tableaux, Gelfand
   patterns, and branching rules for classical groups}. J. Algebra 
{\bf 164}, 299--360 (1994)

\itemitem{[RZ]}
{Retakh, V., Zelevinsky, A.:}
{Base affine space and 
canonical basis in
irreducible representations of $Sp(4)$}.
Dokl. Acad. Nauk USSR {\bf 300}, 31--35  (1988)

\itemitem{[S]}
{Shtepin, V. V.:}
 {Intermediate Lie algebras and their
finite-dimensional representations}. Russian Akad.
Sci. Izv. Math.
{\bf 43}, 559--579  (1994)

\itemitem{[T]}
{Tarasov, V. O.:}
{Irreducible monodromy matrices for the $R$-matrix of the
$XXZ$-model and lattice local quantum Hamiltonians}. {Theor. Math. Phys.}
{\bf 63},
440--454 (1985)

\itemitem{[W]}
{Weyl, H.:}
{Classical Groups, their Invariants and Representations}.
Princeton NJ:
{Princeton Univ. Press}
1946

\itemitem{[Wo]}
{Wong, M. K. F.:}
{Representations of the orthogonal
group. I. Lowering and raising operators of the
orthogonal group and matrix elements of the generators}.
{J. Math. Phys.}
{\bf 8}, 1899--1911  (1967)

\itemitem{[WY]}
{Wong, M. K. F., Yeh, H.-Y.:}
{The most degenerate irreducible representations
of the symplectic group}.
{J. Math. Phys.}
{\bf 21}, 630--635  (1980)

\itemitem{[Z1]} 
{Zhelobenko, D. P.:}
{The classical groups. Spectral
analysis of their
finite-dimensional representations}. Russ. Math. Surv. 
{\bf 17}, 1--94  (1962)

\itemitem{[Z2]}
{Zhelobenko, D. P.:}
{Compact Lie groups and their representations}.
{Transl. of Math. Monographs}
{\bf 40} AMS,
Providence RI
1973

\itemitem{[Z3]} 
{Zhelobenko, D. P.:}
{$S$-algebras and Verma modules
   over reductive Lie algebras}.
Soviet. Math. Dokl.
{\bf 28}, 696--700  (1983)

\itemitem{[Z4]} 
{Zhelobenko, D. P.:}
{$Z$-algebras over reductive Lie
   algebras}. Soviet. Math. Dokl.
{\bf 28}, 777--781  (1983)

\itemitem{[Z5]} 
{Zhelobenko, D. P.:}
{On Gelfand--Zetlin
   bases for classical Lie algebras}. In:
Kirillov,~A.~A. (ed.)
 Representations of Lie groups and
   Lie algebras, pp. 79--106.
Budapest: Akademiai Kiado  1985

\itemitem{[Z6]} 
{Zhelobenko, D. P.:}
{Extremal projectors and
generalized Mickelsson algebras on reductive Lie algebras}.
Math. USSR-Izv. {\bf 33}, 85--100  (1989)

\enddocument